\documentclass[10pt]{amsart}
\theoremstyle{plain}
\usepackage{amsmath, amsfonts, amsthm}
\usepackage{graphics}
\usepackage{color}
\usepackage{epsfig}
\usepackage[all]{xy}
\title{Grid Diagrams and Legendrian Lens Space Links}

\author{Kenneth L.\ Baker}
\author{J.\ Elisenda Grigsby}

\theoremstyle{plain}
\newtheorem{theorem}{Theorem}[section]
\newtheorem{lemma}[theorem]{Lemma}
\newtheorem{proposition}[theorem]{Proposition}
\newtheorem{corollary}[theorem]{Corollary}

\newtheorem*{theoremA}{Theorem~\ref{thm:Legrealize}}
\newtheorem*{corB}{Corollary~\ref{cor:gridnumber}}
\newtheorem*{bergeconj}{Berge Conjecture}
\theoremstyle{definition}
\newtheorem{definition}[theorem]{Definition}
\newtheorem{remark}[theorem]{Remark}

\newcommand{\Ozsvath}{{Ozsv{\'a}th} }

\newcommand{\bdry}{\ensuremath{\partial}}

\newcommand{\Q}{\ensuremath{\mathbb{Q}}}
\newcommand{\R}{\ensuremath{\mathbb{R}}}
\newcommand{\Z}{\ensuremath{\mathbb{Z}}}

\newcommand{\C}{\ensuremath{\mathbb{C}}}

\newcommand{\OO}{\ensuremath{\mathbb{O}}}
\newcommand{\XX}{\ensuremath{\mathbb{X}}}

\newcommand{\st}{\ensuremath{\mbox{\tiny{st}}}}
\newcommand{\UT}{\ensuremath{\mbox{\tiny{UT}}}}
\newcommand{\tb}{\ensuremath{\mbox{\rm t\!b}}}

\newcommand{\GN}{\ensuremath{\mbox{\sc gn}}}

\newcommand{\tagGNOneMatsuda}{$\star$}

\begin{document}
\bibliographystyle{halpha}

\begin{abstract} 
Grid diagrams encode useful geometric information about knots in $S^3$.  In particular, they can be used to combinatorially define the knot Floer homology of a knot $K \subset S^3$ \cite{GT0607691,GT0610559}, and they have a straightforward connection to Legendrian representatives of $K \subset (S^3, \xi_{st})$, where $\xi_{st}$ is the standard, tight contact structure \cite{MR2240685,GT0611841}.  The definition of a grid diagram was extended, in \cite{GT07100359}, to include a description for links in all lens spaces, resulting in a combinatorial description of the knot Floer homology of a knot $K \subset L(p,q)$ for all $p \neq 0$.  In the present article, we explore the connection between lens space grid diagrams and the contact topology of a lens space.  Our hope is that an understanding of grid diagrams from this point of view will lead to new approaches to the Berge conjecture, which claims to classify all knots in $S^3$ upon which surgery yields a lens space.
\end{abstract}
\maketitle

\section{Introduction} 
\begin{figure}
\begin{center}
\input{trefoilgriddiagram.pstex_t}\\
\vskip.5cm
\includegraphics[height=2in]{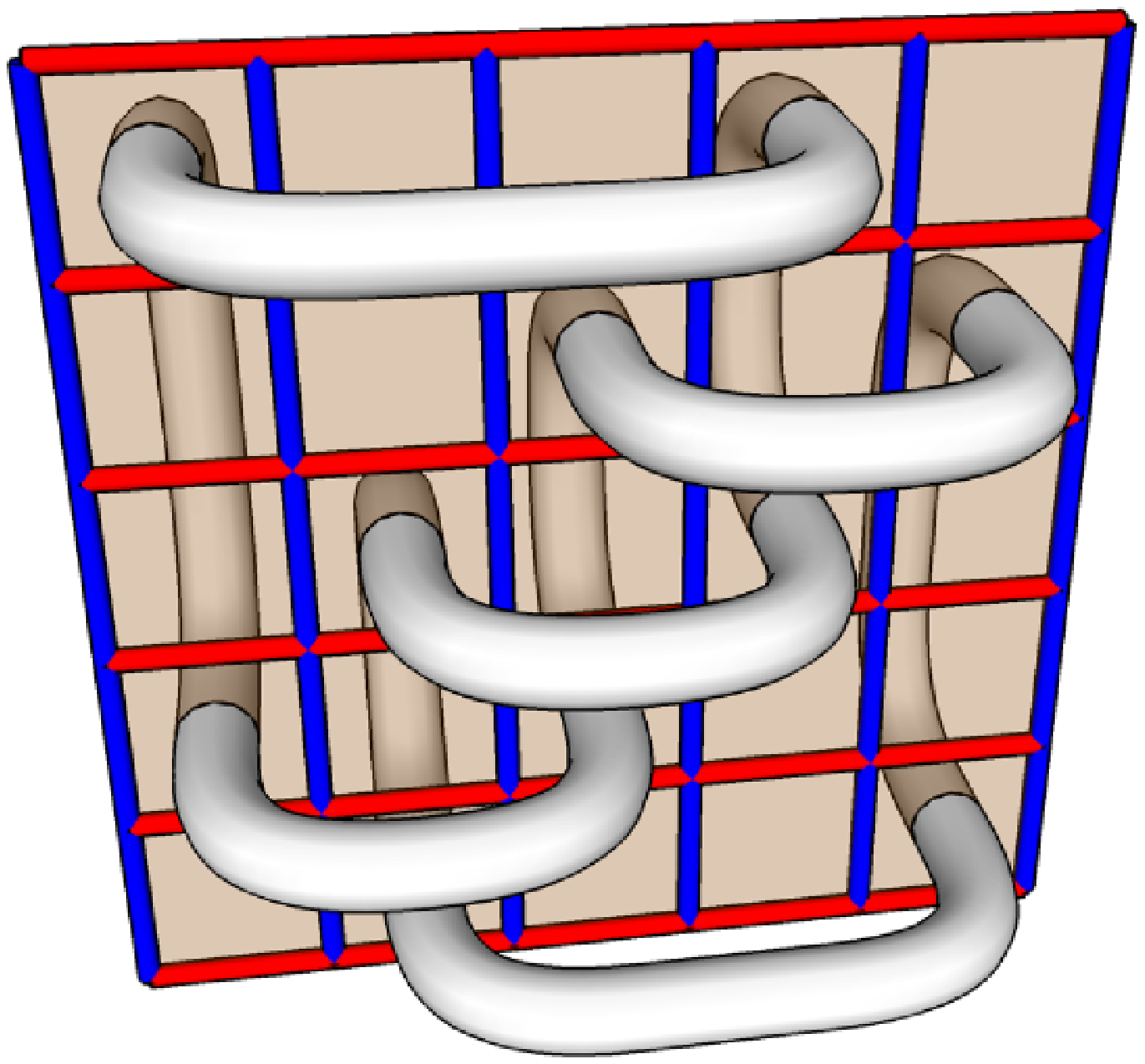} \quad\quad\quad \quad \includegraphics[height=2in]{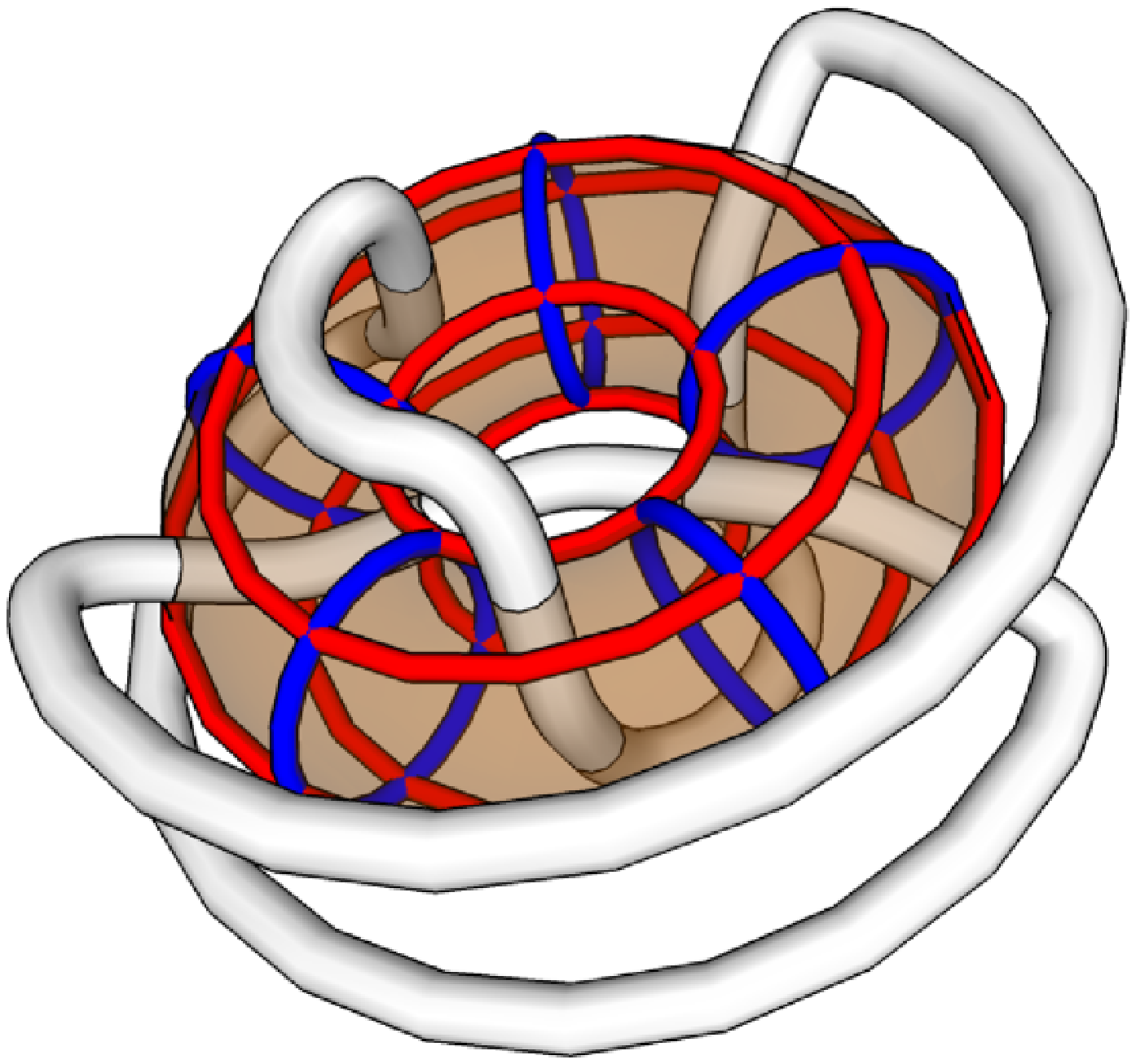} 
\end{center} 
\caption{An $S^3$ grid diagram for a trefoil, along with two views of the knot it represents.}
\label{fig:trefoilgrid}
\end{figure}

A grid diagram provides a simple combinatorial means of encoding the data of a link in $S^3$, as in Figure~\ref{fig:trefoilgrid}.  Though grid diagrams first made an appearance in the late 19th century \cite{Brunn}, they have enjoyed an abundance of recent attention, due primarily to their connection to contact topology \cite{MR555835,MR1339757,MR2232855,MR2240685,GT0611841} and combinatorial Heegaard Floer homology \cite{GT0607691, GT0610559}.

The definition was extended in \cite{GT07100359} to provide a means of encoding  the data of all lens space links via grid diagrams, leading to a combinatorial description of the knot Floer homology of a lens space knot.  Figure~\ref{fig:twistedgrid} illustrates the notion of a lens space grid diagram; we delay its precise definition until Section~\ref{sec:LegendReps}.

In the present article, we explore the connection between the contact geometry of a lens space and grid diagrams for lens space links.  
\begin{theoremA} 
Associated to a grid diagram for a link $K$ in the lens space $L(p,q)$ is a unique Legendrian representative of $m(K)$ with respect to $(-L(p,q),\xi_{\UT})$.
\end{theoremA}
Here, $m(K)$ denotes the topological mirror of $K$ and $\xi_{UT}$ is a canonical co-oriented contact structure on $-L(p,q)$. See Section \ref{sec:orientation} for a detailed discussion of notation and orientation conventions.

  Our approach is similar in spirit to that taken by Matsuda and Menasco in \cite{GT07082406}.   In a very natural geometric fashion, a grid diagram describing a lens space link corresponds to a Legendrian representative of the link with respect to a canonical cooriented tight contact structure on the lens space.  A contact structure $\xi$ on a $3$--manifold $M$ is a smooth, nowhere integrable $2$--plane field.  A link in the contact manifold $(M,\xi)$ is {\em Legendrian} if it is everywhere tangent to $\xi$ and {\em transverse} if it is everywhere transverse to $\xi$.  A contact structure $\xi$ is said to be {\em tight} if the linking number of any trivial Legendrian knot with its contact push-off is negative (equivalently, $(M,\xi)$ contains no overtwisted disks).  Furthermore, a tight contact structure $\xi$ on a $3$--manifold $M$ is {\em universally tight} if the lift of $\xi$ to the universal cover of $M$ is also tight.

By Honda's classification \cite{MR1786111}, each lens space $L(p,q)$ has two distinct (positive, co-oriented) universally tight contact structures when $0 < q < p-1$ and just one when $q=p-1$.  Given one universally tight contact structure on $L(p,q)$, the other is obtained by flipping its co-orientation.  When $q=p-1$, these two are isotopic.

Suppressing the specific lens space from the notation, we will denote a particular co-oriented universally tight contact structure on $L(p,q)$ by $\xi_{\UT}$.  We explicitly construct $\xi_{\UT}$ in Section~\ref{sec:universallytight}.

We use {\em toroidal fronts} to mediate between grid diagrams and Legendrian lens space links.  Just as planar fronts uniquely specify Legendrian links in $(\R^3, \xi_{st})$ via projection to the $xz$--plane, toroidal fronts uniquely specify Legendrian links in $(L(p,q),\xi_{\UT})$ via projection to a standard Heegaard torus.  After defining toroidal fronts and their correspondence with Legendrian lens space links in Section~\ref{sec:toroidalfront}, we detail the relationship between toroidal fronts and grid diagrams in Sections~\ref{sec:LegendReps} and \ref{sec:Equiv}.

\subsection{Motivation}
One motivation for developing the connection between grid diagrams and Legendrian lens space links (and their toroidal fronts) is to provide the foundation for a systematic study of Legendrian links in a collection of contact manifolds other than $(S^3,\xi_{\st})$.  Furthermore, certain elements in the knot Floer homology chain complex of (the mirror of) a lens space knot associated to a grid diagram should yield powerful invariants of the corresponding Legendrian or transverse representatives of the knot.  
In particular, one ought to be able to use such invariants to detect transversely non-simple knots in lens spaces, as Ng, Ozsv{\'a}th, and Thurston did for knots in $S^3$, \cite{GT0703446}.

Moreover, this work reveals an inlet for the use of contact structures in studying a fundamental question about the interconnectedness of $3$--manifolds through Dehn surgery on knots.  Lens spaces are precisely the manifolds that may be obtained by Dehn surgery on the unknot in $S^3$.  The resolutions of The Knot Complement Problem \cite{gl:kadbtc} and Property R \cite{gabai:fatto3m} show that the unknot is the only knot in $S^3$ on which a non-trivial Dehn surgery may produce the lens spaces $S^3$ and $S^1 \times S^2$, respectively.    For many other lens spaces, this is not the case.  Indeed, all torus knots admit lens space surgeries \cite{moser}, as do certain cables of torus knots (the only satellite knots admitting lens space surgeries) \cite{bleilerlitherland} and many hyperbolic knots \cite{bailyrolfsen}, \cite{ronandron}, \cite{Berge}.  The Berge Conjecture proposes a classification of all knots in $S^3$ and their Dehn surgeries that yield lens spaces.  Originally stated in terms of a homotopy condition for a knot embedded on the surface of a genus $2$ Heegaard surface in $S^3$, the Berge Conjecture is most succinctly stated as follows:

\begin{bergeconj}[\cite{Berge}] 
If a knot $K$ in a lens space $L(p,q)$, with $p  \not\in \{0, 1\}$, admits an $S^3$ Dehn surgery, then it has grid number $1$.   
\end{bergeconj}

The {\em grid number}, \GN, of a link is the minimum grid number over all grid diagrams representing the link.  Figure~\ref{fig:twistedgrid} clarifies the definition of the grid number of a grid diagram.  In light of the above, one might hope for an upper bound on {\sc gn} for knots admitting $S^3$ surgeries.

For links in $S^3$, Matsuda \cite{MR2240685} proves 
\[\GN(K) \geq -\overline{\tb}(K) - \overline{\tb}(m(K)).\]
Here, $\tb$ denotes the classical Thurston-Bennequin number associated to a Legendrian link in $(S^3,\xi_{\st})$, and $\overline{\tb}(K)$ denotes the maximal Thurston-Bennequin number over all Legendrian representatives of $K$.  

Though this provides a lower bound on grid number, Ng speculates, in \cite{GT0612356}, that this bound may be sharp.  A proof of this would imply that $\overline{\tb}(K)$ and $\overline{\tb}(m(K))$ determine $\GN(K)$, and thus an upper bound on the maximal Thurston-Bennequin numbers of a link and its mirror would produce an upper bound on the grid number.

For links in lens spaces, we prove an analogue of Matsuda's bound in Section~\ref{sec:ClassInv}.  
\begin{corB}
\[\GN(K) \geq -\overline{\tb}(K) - \overline{\tb}(m(K))\]
for each link $K \subset L(p,q)$.
\end{corB}
\noindent This requires an extension of the definition of $\tb(K)$ to all Legendrian links in contact rational homology spheres, provided in Definition~\ref{defn:tb}.

 If Matsuda's bound is sharp for lens space links, one can bound from above the grid number of a knot admitting an $S^3$ surgery by finding upper bounds for the Legendrian contact invariants $\overline{\tb}(K)$ and $\overline{\tb}(m(K))$.  The hope is that suitably understanding Legendrian representatives of knots in $L(p,q)$ with $S^3$ surgeries will shed light on the Berge Conjecture.

\subsection{Acknowledgments}
We thank John Baldwin, John Etnyre, Matt Hedden, Lenny Ng, Peter \Ozsvath, and Andras Stipsicz for many enlightening conversations.  We are especially grateful to Andras Stipsicz, who provided a key observation crucial in the proof of Corollary \ref{cor:zeroorone}.  The first author was partially supported by NSF Grant DMSÐ 0239600; the second author was partially supported by an NSF postdoctoral fellowship.  

\section{Notation and Orientation Conventions} \label{sec:orientation}
Throughout the paper, $L(p,q)$ will denote $-\frac{p}{q}$ surgery on the unknot, where $p$ and $q$ are coprime integers such that $0\leq |q| < p$.  We view $S^3$ as the lens space $L(1,0)$ and will not consider $S^1 \times S^2$.  

Let $G_K$ be a grid diagram representing a link $K \subset L(p,q')$.  Then $G_K$ will naturally yield a Legendrian representative, which we will denote $\mathcal{L}_K$, of the topological mirror of $K$, $m(K) \subset (L(p,q), \xi_{UT})$.   Here, $q$ satisfies $qq' \equiv -1 \mod p$.  Note that, by the classification of lens spaces up to orientation-preserving diffeomorphism, $L(p,q) \cong -L(p,q')$.  We describe $\xi_{UT}$ as the kernel of a globally-defined $1$--form in the next section.

The correspondence between a grid diagram (or, more generally, a Heegaard diagram) compatible with a knot, $K$, and a Legendrian representative of $m(K)$, though odd, is by now standard in the literature.  See, e.g., \cite{GT0611841}, where a grid diagram for $K \subset S^3$ is associated to a Legendrian representative of $K$ with respect to $\xi_{\st}$ on $S^3$ with the opposite orientation (i.e., a Legendrian representative of the mirror of $K$).  See also \cite{MR2153455}, which defines the contact invariant of a fibered link in $M$ as an element of the Heegaard Floer homology of $-M$.

This will have the unfortunate effect that the orientation on a Heegaard torus associated to a grid diagram is opposite to the orientation on a Heegaard torus associated to a toroidal front diagram.  This situation, though confusing, is unavoidable, since it is important to match the existing convention in the literature.  To avert confusion, it will be convenient to define the notion of a {\em dual grid diagram} (Definition \ref{defn:dualgriddiagram}), a grid diagram for the mirror of a given link, $K$.  In fact, after briefly recalling the original definition of a grid diagram of a link, we will thereafter work exclusively with dual grid diagrams---which may be canonically identified with toroidal front diagrams---throughout Section \ref{sec:LegendReps}.  We return to working with the original grid diagrams in Section \ref{sec:Equiv}, after we have proven a correspondence between planar subsets of toroidal fronts and planar fronts in Section \ref{sec:planarfront}.  The coordinates on a dual grid diagram will always be the $(\theta_1,\theta_2)$ coordinates inherited from the quotient map $S^3 \to L(p,q)$ described in the next section.

\section{$\xi_{UT}$ and Toroidal Front Diagrams for Lens Spaces}
\subsection{Construction of $\xi_{UT}$} \label{sec:universallytight}
We begin with a construction of the universally tight contact structure $\xi_{\UT}$ on $L(p,q)$.   Whenever we refer to the contact structure $\xi_{\UT}$ on $L(p,q)$ we will mean the isotopy class of the one constructed as follows.  

Thinking of $S^3$ as the unit sphere in $\C^2$,
\[S^3 = \left\{(u_1,u_2) \in \C^2 = \C_1 \times \C_2| \,\, |u_1|^2 + |u_2|^2 = 1\right\},\] 
the standard tight contact structure on $S^3$ is given by (cf.\ \cite{MR2194671}) 
\[\xi_{\st} = \ker \alpha,\]
 where $\alpha$ is the $1$-form 
 \[\alpha = r_1^2 \, d\theta_1 +r_2^2 \, d\theta_2\] in terms of polar coordinates $(r_1,\theta_1)= u_1 \in \C_1$, $(r_2,\theta_2)=u_2 \in \C_2$.
Thinking of $S^3$ as $\R^3 \cup \{\infty\}$, it is natural to identify the circle corresponding to $r_1 = 1$ with the $z$--axis $\cup \{\infty\}$ and the circle corresponding to $r_2 = 1$ with the unit circle on the $xy$--plane.

Let $\omega_p = e^{\frac{2\pi i}{p}}$.  Then $L(p,q)$ can be identified as the quotient:
\[L(p,q) = S^3/(u_1,u_2)\sim(\omega_p u_1, \omega_p^{q}u_2).\] 
 Noting that 
 \[\left\{(r_1,\theta_1,r_2,\theta_2)\,\,|\,\, \theta_1 \in [0,2\pi),\theta_2 \in \left[0,\frac{2\pi}{p}\right)\right\} \subset S^3\]
 is a fundamental domain for the $\mathbb{Z}_p$ action of $\sim$, and that the coordinate $r_2$ may be recovered from the condition $r_1^2 + r_2^2 = 1$, we will specify points in $L(p,q)$ by: 
 \[\left\{(r_1,\theta_1,\theta_2)\,\,|\,\,r_1 \in [0,1],\theta_1 \in [0,2\pi), \theta_2 \in \left[0,\frac{2\pi}{p}\right)\right\}.\]

Let $\pi\colon S^3 \to L(p,q)$ be the covering map induced by the $\sim$ equivalence.
 Then $\xi_{\UT}=\pi_*(\xi_{\st})$ is a well-defined tight contact structure on $L(p,q)$, since the $1$-form $\alpha$ is constant along tori of constant radius $r_1$ and these tori are mapped to themselves under the $\mathbb{Z}_p$ action of $\sim$.  
 
 Note that the global $1$--form $\alpha$ induces a co-orientation on $(S^3,\xi_{\st})$ and hence on $(L(p,q),\xi_{UT})$.  The other universally tight contact structure on $L(p,q)$ (for $0< q < p-1$) may be obtained by using $-\alpha$ in the above construction.
Regardless, most of this paper is insensitive to the co-orientation.

\subsection{Toroidal Front Diagrams} \label{sec:toroidalfront}
Recall that a smooth link in a contact $3$--manifold $(M,\xi)$ is said to be {\em Legendrian} if its tangent vectors are everywhere tangent to the contact planes.  A Legendrian isotopy is a smooth isotopy through Legendrian links, and the terminology {\em Legendrian link} refers to the Legendrian isotopy class of a Legendrian link as well as a specific representative.

The radial projection of a Legendrian link in $(L(p,q),\xi_{\UT})$ onto the radius $r_1 =  \frac{1}{\sqrt{2}}$ Heegaard torus
\[\Sigma = \left\{(r_1,\theta_1,\theta_2)| \,\, r_1 = \frac{1}{\sqrt{2}}\right\} \subset L(p,q)\]
gives a toroidal front diagram (defined below) on $\Sigma$ from which we may recover the Legendrian link.  Observe that $\Sigma$ separates $L(p,q)$ into the two solid tori  
\[
V^\alpha = \left\{(r_1,\theta_1,\theta_2)\,\,|\,\, r_1 \in \left[0, \frac{1}{\sqrt{2}}\right]\right\} \mbox{ and }
V^\beta = \left\{(r_1,\theta_1,\theta_2)\,\,|\,\, r_1 \in \left[\frac{1}{\sqrt{2}}, 1\right]\right\}
\]
which are oriented so that $\Sigma = \bdry V^\alpha = -\bdry V^\beta$.  (In terms of the standard surgery description of $L(p,q)$, we may view $V^\alpha$ as the surgered neighborhood and $V^\beta$ as the exterior of the unknot in $S^3$.)  

\begin{remark}  We use the notation $\alpha$ both to refer to the $1$--form defining $\xi_{\st}$ and to describe various objects in one of the two solid tori bounded by $\Sigma$.  These two usages are common in the literature, and context should prevent confusion.
\end{remark}

Points on $\Sigma$ may be uniquely specified by the fundamental domain 
\[\left\{(\theta_1,\theta_2)\,\,|\,\,\theta_1 \in [0,2\pi), \theta_2 \in \left[0,\frac{2\pi}{p}\right)\right\},\]
the intersection of a fundamental domain for $L(p,q)$ with $\Sigma$.
Then with respect to the bases
\[\left\{\frac{\partial}{\partial \theta_1}, \frac{\partial}{\partial \theta_2}\right\} \mbox{ and } \{d\theta_1,d\theta_2\}\]
for the tangent and cotangent space at each point,
$\Sigma$ is naturally parallelizable.

We are now ready to define toroidal front diagrams.

\begin{definition}
A {\em toroidal front diagram} (or simply a {\em toroidal front}) for a Legendrian link $\mathcal{L}$ in $(L(p,q), \xi_{\UT})$ is an immersion $f: S^1 \amalg \ldots \amalg S^1 \rightarrow \Sigma$ with the following properties:
\begin{itemize}
\item $f$ is an embedding except at finitely many transverse double points and smooth except at finitely many cusps.
\item The slopes $\frac{d\theta_2}{d\theta_1}$ of the tangent vectors at the smooth points satisfy $\frac{d\theta_2}{d\theta_1} \in (-\infty,0)$.
\item Each cusp is semi-cubical.  See the discussion following Lemma 2.45 in \cite{MR2194671} for the precise definition of a semi-cubical cusp.
\end{itemize}
\end{definition}

\begin{proposition}\label{prop:frontgivesleg}
A toroidal front diagram uniquely specifies a Legendrian link up to Legendrian isotopy.
\end{proposition}
\begin{proof}
Let $\gamma$ be the image of $f$ on $\Sigma$.  We claim that $\gamma$ is naturally the $(\theta_1,\theta_2)$ projection of a Legendrian link in $(L(p,q),\xi_{\UT})$, where the $r_1$ coordinate is recovered from the Legendrian condition.  

Set $m = \frac{d\theta_2}{d\theta_1}$.  Then the condition that the vectors tangent to a Legendrian curve lie in $\ker \alpha$ implies that $r_1^2 d\theta_1 + (1-r_1^2) d\theta_2=0$ and hence $\frac{r_1^2}{(r_1^2-1)} = m$.  This gives the unique  non-negative solution $r_1 = \sqrt{\frac{m}{m-1}}$.  Because $m \in (-\infty,0)$ at the smooth points of $\gamma$, we obtain $r_1 \in (0,1)$ there.  At the cusps, the one-sided tangencies agree, giving a slope $m \in (-\infty,0)$ and hence a radius $r_1 \in (0,1)$ there, too.  

The condition that $f$ is smooth away from the cusps ensures that the corresponding Legendrian link is smooth away from the preimages of the cusps.  The condition that the cusps of $f$ are semi-cubical ensures that the Legendrian link is smooth in a neighborhood of the preimage of the cusps.  
\end{proof}

\begin{proposition}\label{prop:leggivesfront}
Every Legendrian isotopy class of Legendrian links in $(L(p,q),\xi_{\UT})$ has a representative admitting an associated toroidal front projection.
\end{proposition}

\begin{proof}
The cores of the Heegaard solid tori $V^\alpha$ and $V^\beta$ correspond to the circles where $r_1=0$ and $r_1=1$, respectively.  By a Legendrian isotopy, a Legendrian link may be made disjoint from these two cores.  This ensures that the Legendrian link has a well-defined projection $\Pi \colon (r_1,\theta_1,\theta_2) \mapsto (\theta_1,\theta_2)$ to the Heegaard torus $\Sigma$.  By a further Legendrian isotopy,
we may ensure that the image of the projection is an embedding except at finitely many transverse double points.  Such Legendrian isotopies exist because the subspace of Legendrian representatives in a particular Legendrian isotopy class disjoint from the two cores and having a generic projection, as above, is open, dense, and positive-dimensional.  
 We claim that this projection is a toroidal front diagram for the Legendrian link.

Because each component of the Legendrian link is a smooth curve, its projection under $\Pi$ is smooth except for where the link has a tangent vector that is parallel to $\frac{\partial}{\partial r_1}$.  For the points where the projection is smooth, the Legendrian condition $r_1^2 d\theta_1 + (1-r_1)^2 d\theta_2=0$ implies that $\frac{d\theta_2}{d\theta_1} \in (-\infty,0)$ since $r_1 \in (0,1)$.  For the points where the projection is not smooth, the Legendrian condition $r_1^2 d\theta_1 + (1-r_1)^2 d\theta_2=0$ implies that there are local coordinates for the projection presenting a neighborhood of the non-smooth point as a semi-cubical cusp \cite{MR2194671}.  Since the projection is compact, there can only be finitely many cusps.
\end{proof}

\begin{remark}
Note that toroidal fronts for Legendrian links in lens spaces behave locally much like planar fronts for Legendrian links in $\R^3$.  For example: If two arcs $a_1$ and $a_2$ of a toroidal front transversally intersect at a point $P$ with slopes $m_1$ and $m_2$ respectively such that $-\infty<m_1 < m_2<0$, then the point projecting to $P$ on $a_1$ lies in front of (has greater $r_1$ than) the point projecting to $P$ on $a_2$.  The correspondence between planar subsets of toroidal front projections (which represent Legendrian tangles in $(\R^3,\xi_{\st})$) and standard planar front projections for Legendrian tangles in $(\R^3,\xi_{\st})$ is made explicit in Section \ref{sec:planarfront}.
\end{remark}

\section{Grid Diagrams and Legendrian Knots}\label{sec:LegendReps}
Let us quickly remind the reader of the definition of a toroidal grid diagram $G_K$ for a link $K$ in $L(p,q')$.

\begin{figure}
\begin{center}
\input{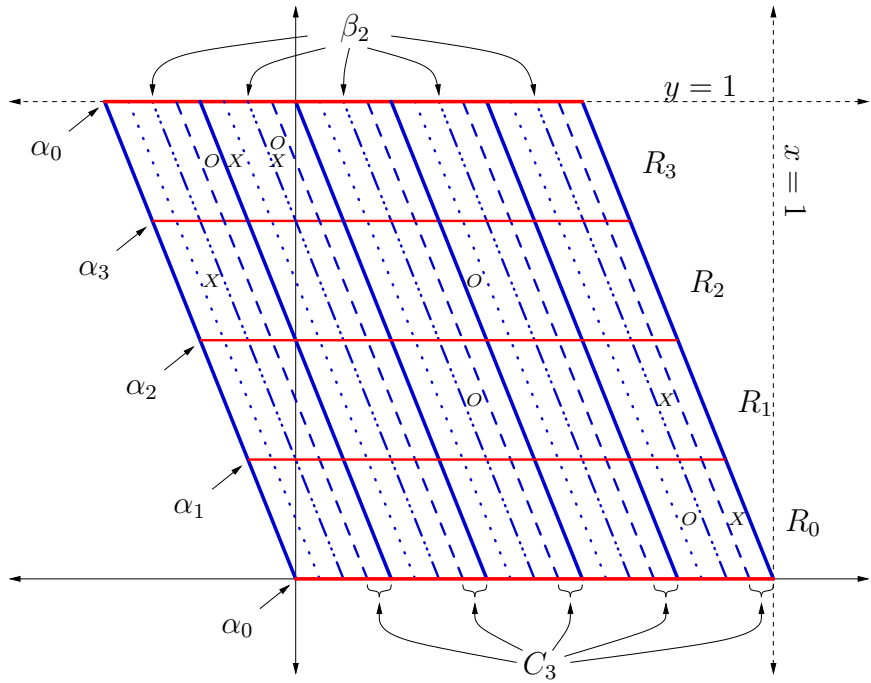}
\end{center}
\caption{An example of a toroidal grid diagram $G_K$ with grid number $n=4$ for a link $L$ in $L(5,2)$.  Here, $C_3$ is one of the four columns, while $R_i$ are the rows.}
\label{fig:twistedgrid}
\end{figure}

\begin{definition} \label{def:TTGridDiag}  A {\em (twisted toroidal) grid diagram} $G_K$ with grid number $n$ for $L(p,q')$ consists of  a five-tuple $(T^2,\vec{\alpha},\vec{\beta},\vec{\mathbb{O}},\vec{\mathbb{X}})$, illustrated in Figure~\ref{fig:twistedgrid}, where: 
  \begin{itemize}
	  \item $T^2$ is the standard oriented torus $\R^2 / \Z^2$, identified with the quotient of $\mathbb{R}^2$ (with its standard orientation) by the $\mathbb{Z}^2$ lattice generated by the vectors $(1,0)$ and $(0,1)$. \vspace{2mm}
    \item $\vec{\alpha} = \{\alpha_0, \ldots, \alpha_{n-1}\}$ are the $n$ images $\alpha_i$ in $T^2 = \mathbb{R}^2/\mathbb{Z}^2$ of the lines $y = \frac{i}{n}$ for $i \in \{0, \ldots n-1\}$.  Their complement $T^2-\alpha_0 - \ldots - \alpha_{n-1}$ has $n$ connected annular components, which we call the {\em rows} of the grid diagram.  \vspace{2mm}
    \item $\vec{\beta} = \{\beta_0, \ldots, \beta_{n-1}\}$ are the $n$ images $\beta_i$ in $T^2 = \mathbb{R}^2/\mathbb{Z}^2$ of the lines $y = -\frac{p}{q'}(x-\frac{i}{pn})$ for $i \in \{0, \ldots n-1\}$.  Their complement $T^2 - \beta_0 - \ldots - \beta_{n-1}$ has $n$ connected annular components, which we call the {\em columns} of the grid diagram.  \vspace{2mm} 
    \item $\vec{\mathbb{O}} = \{O_0, \ldots, O_{n-1}\}$ are $n$ points in $T^2 - \vec{\alpha} - \vec{\beta}$ with the property that no two $O$'s lie in the same row or column.  \vspace{2mm}
    \item $\vec{\mathbb{X}} = \{X_0, \ldots, X_{n-1}\}$ are $n$ points in $T^2 - \vec{\alpha} - \vec{\beta}$ with the property that no two $X$'s lie in the same row or column. \vspace{2mm}
  \end{itemize}
\end{definition}

We refer to the connected components of $T^2 - \vec{\alpha} - \vec{\beta}$ as the {\em fundamental  parallelograms} of the grid diagram.

Two five-tuples $(T^2, \vec{\alpha},\vec{\beta},\vec{\OO},\vec{\XX})_1$ and $(T^2,\vec{\alpha},\vec{\beta},\vec{\OO},\vec{\XX})_2$ are equivalent (and, hence, represent the same grid diagram $G_K$) if there exists an orientation-preserving diffeomorphism $(T^2)_1 \rightarrow (T^2)_2$ respecting the markings (up to cyclic permutation of their labels).

One associates a unique oriented link $K$ in $L(p,q')$ to a grid diagram $G_K=(T^2, \vec{\alpha},\vec{\beta},\vec{\OO},\vec{\XX})$ as follows: 
\begin{enumerate}
\item First attach solid tori $V^\alpha$ and $V^\beta$ to the torus $T^2$ so that $T^2 = \bdry V^\alpha = -\bdry V^\beta$ the $\alpha$ curves of $G_K$ are meridians of $V^\alpha$ and the $\beta$ curves of $G_K$ are meridians of $V^\beta$.  This forms the standard embedding of $T^2$ with its decorations into $L(p,q')$.
\item Connect each $X_i$ to the unique $O_j$ lying in the same row as $X_i$ by an oriented ``horizontal'' arc embedded in that row of $T^2$, disjoint from the $\vec{\alpha}$ curves.
\item Next connect each $O_{j}$ to the unique $X_m$ lying in the same column as $O_j$ by an oriented ``slanted'' arc embedded in that column of $T^2$, disjoint from the $\vec{\beta}$ curves.\footnote{If an $O_i$ and an $X_j$ coincide, then we may take the slanted arc connecting them to be trivial joining with the horizontal arc to form a full circle.} 
\item The union of these two collections of $n$ arcs forms an immersed (multi)curve $\gamma$ in $T^2$.  Remove all self-intersections of $\gamma$ by pushing the interiors of the horizontal arcs slightly down into $V^\alpha$ and the interiors of the slanted arcs slightly up into $V^\beta$.
\end{enumerate}

\begin{figure}
\centering
\input{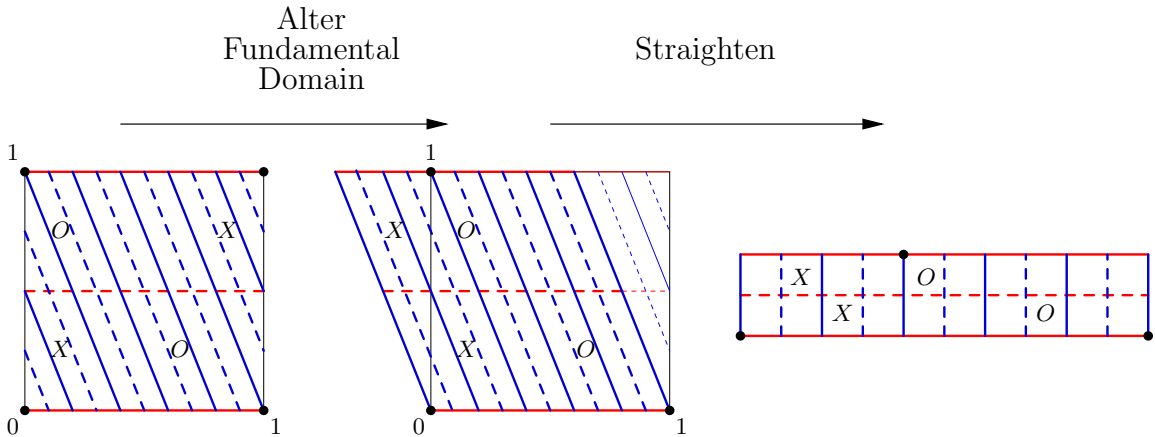}
\caption{A grid number $2$ diagram for a link in $L(5,2)$ is shown before and after straightening.  One obtains the torus represented by the grid diagram $G_K$ by identifying opposite edges of the square on the left-hand side in the standard manner.  One can construct the same object, beginning with the rectangle on the right-hand side, by first identifying the two vertical edges and then the two boundary circles of the resulting annulus with an appropriate twist so that the thick black dots all represent the same point after the identifications.}
\label{fig:StraightT2}
\end{figure}

It will often be convenient to pick a particular fundamental domain for the grid diagram and ``straighten'' it out so that the image of the $\alpha$ curves are horizontal, the $\beta$ curves are vertical, and each row is connected in the fundamental domain, as in Figure \ref{fig:StraightT2}.  From now on, we will always represent a grid diagram in this manner.

In accordance with Section~\ref{sec:orientation}, we make the following definition in order to identify a grid diagram for a link $K$ in $L(p,q')$ with a Legendrian representative of the link in the contact manifold $(L(p,q),\xi_{\UT})$.  As before, $qq' \equiv -1 \mod p$.

\begin{definition} \label{defn:dualgriddiagram} 
Given a grid diagram $G_K$ representing a link in $L(p,q')$, let $G_K^*$ denote the {\em dual grid diagram}, obtained as follows and illustrated in Figure~\ref{fig:DualGrid}.
  \begin{enumerate}
    \item Begin with any straightened fundamental domain for $G_K$.
    \item Rotate the fundamental domain $90^{\circ}$ clockwise; the $\alpha$  curves of the old fundamental domain are now vertical, and the $\beta$ curves are now horizontal.
    \item Chop the rotated fundamental domain along the newly horizontal arcs of the bottom-most $\beta$ curve.  Then reglue the resulting pieces together appropriately along the newly vertical arcs of the $\alpha$ curve on their sides.  This produces a new straightened fundamental domain for the torus whose rows are connected. 
    \item Relabel all horizontal circles on the torus as  $\alpha$ circles and all vertical circles as $\beta$ circles, and vice versa.
    \item Relabel all $X$'s as $O$'s and vice versa.
    \item The decorated, straightened fundamental domain represents $G_K^*$.
  \end{enumerate}
\end{definition}

\begin{figure}
\centering
\input{DualGrid.pstex_t}
\caption{}
\label{fig:DualGrid}
\end{figure}

Note that $G_K^*$ defines a link in $L(p,q)$, by the same procedure described above (connect $X$'s to $O$'s in each row and $O$'s to $X$'s in each column, with vertical passing over horizontal).

\subsection{Identifying $G_K^*$ with the Constant Radius Heegaard Torus $\Sigma$.}\label{sec:TonSigma}
The torus arising in the definition of a grid diagram is, in fact, a Heegaard torus associated to a Heegaard decomposition of the appropriate lens space.  Therefore, if we begin with a grid diagram $G_K$ representing a link $K \subset L(p,q')$, then it is natural to identify the torus associated to the dual grid diagram, $G_K^*$, with the torus 
\[\Sigma = \left\{(r_1,\theta_1,\theta_2) \,\,|\,\,r_1 = \frac{1}{\sqrt{2}}, \theta_1 \in [0,2\pi), \theta_2 \in \left[0,\frac{2\pi}{p}\right)\right\} \subset L(p,q)\]
defined in Section~\ref{sec:toroidalfront} (where the coordinates give an implicit choice of fundamental domain for $\Sigma$) 
so that each horizontal $\alpha$ curve of slope $0$ corresponds to a circle of constant $\theta_2$ and that each slanted $\beta$ curve of slope $-\frac{p}{q}$ 
corresponds to a circle of constant $\theta_1$ taken mod $\frac{2\pi}{p}$.  (In the fundamental domain the $\beta$ curves correspond to a union of the $p$ lines in $\Sigma$ with $\theta_1$ coordinate in the set ${\bf a} = \left\{a + \frac{2\pi k}{p}\right\}$, for $k \in \{0,\ldots p-1\}$ and some fixed $a \in \left[0,\frac{2\pi}{p}\right)$.)  In particular, the $\alpha$ curves are meridians of $V^\alpha$ and the $\beta$ curves are meridians of $V^\beta$.  

We choose the identification of $T^2$ with $\Sigma$ so that, within each row, the unique $O$ and $X$ have the same $\theta_2$ coordinate and, within each column, the unique $O$ and $X$ have the same $\theta_1$ coordinate mod $\frac{2\pi}{p}$.  
Furthermore, if a pair of an $O$ and an $X$ lie in the same fundamental parallelogram, then it is natural to associate to the pair a trivial Legendrian unknot.\footnote{More precisely, one isotopes the component to lie in a Darboux ball and chooses the Legendrian isotopy class contactomorphic to the $\tb = -1$ unknot in $(S^3, \xi_{\st})$.} In what follows, we may therefore assume, without loss of generality, that no fundamental parallelogram contains both an $X$ and an $O$. 

\begin{remark}
It is worthwhile to remark at this point that there is a natural correspondence between Legendrian links in $(L(p,q), \xi_{\UT})$ and $\Z_p$-symmetric Legendrian links in $(S^3,\xi_{\st})$ obtained via the covering operation.  This correspondence matches the correspondence between grid diagrams $G_K$ for $K \subset L(p,q)$ and their lifts to grid diagrams $G_{\widetilde{K}}$ for $\widetilde{K} \subset S^3$.  See \cite{GT07100359}.
\end{remark}

We are now ready to state and prove the main theorem.  Statements and proofs of many of the supporting lemmas and propositions occupy the remainder of this section:

\begin{theorem} \label{thm:Legrealize} Associated to a grid diagram for a link $K$ in the lens space $L(p,q')$ is a unique Legendrian representative $\mathcal{L}_K$ of $m(K)$ with respect to $(-L(p,q'), \xi_{\UT})$.  Conversely, every Legendrian link $\mathcal{L}_K$ in $(-L(p,q'),\xi_{\UT})$ can be specified by means of a grid diagram for $K$ in $L(p,q')$.
\end{theorem}

\begin{proof}
Let $q$ satisfy $qq' \equiv -1 \mod p$ (see Section \ref{sec:orientation} for a discussion of orientation conventions).  Beginning with $G_K$, a grid diagram for $K \subset L(p,q')$, we produce the dual grid diagram $G_K^*$ associated to $m(K) \subset L(p,q)$ using the procedure described in Definition \ref{defn:dualgriddiagram}.  Lemma \ref{lem:griddiagramtotoroidalfront} then explains how to obtain a toroidal front from a rectilinear projection associated to $G_K^*$, and Proposition \ref{prop:frontgivesleg} associates a unique Legendrian link in $(L(p,q),\xi_{\UT})$ to this toroidal front on $\Sigma$, a coordinatized Heegaard torus for $L(p,q)$.  Proposition \ref{prop:rectilinearchoice} then proves that the choice of rectilinear projection does not affect the Legendrian isotopy class of the resulting Legendrian link.  

Conversely, if we begin with a Legendrian link $\mathcal{L}_K$ representing $m(K) \subset L(p,q)$, Proposition \ref{prop:leggivesfront} associates to it a toroidal front on $\Sigma$.  Lemma \ref{lem:toroidalfronttogriddiagram} then explains how to obtain a grid diagram $G_K^*$ representing $m(K)$ in $L(p,q)$ from the toroidal front.  By reversing the procedure described in Definition \ref{defn:dualgriddiagram}, one obtains a grid diagram $G_K$ representing $K$ in $L(p,q')$.
\end{proof}

\subsection{Constructing a Toroidal Front from a Grid Diagram}\label{subsec:smoothing}
To a grid number $n$ grid diagram we can associate $2^{2n}$ possible piecewise linear projections to $T^2$ (there are $2$ choices for each of the $2n$ horizontal and vertical arcs).  We will call each such projection a {\it rectilinear projection}.

\begin{lemma}\label{lem:griddiagramtotoroidalfront}
A rectilinear projection associated to a grid diagram $G_K^*$ for a link $K \subset L(p,q)$ uniquely specifies a toroidal front for $\mathcal{L}_K$, a Legendrian link in $(L(p,q),\xi_{\UT})$.
\end{lemma}

\begin{proof}
We continue to view the grid diagram $G_K^*$ on $\Sigma$ as described in Section~\ref{sec:TonSigma}.   Recall that we have chosen the identification so that the $O$ and $X$ in each row have the same $\theta_2$ coordinate and the $O$ and $X$ in each column have the same $\theta_1$ coordinate mod $\frac{2\pi}{p}$.  

To obtain a diagram for the link $K$ on $\Sigma$ corresponding to $G_K^*$, we may join each $O$ and $X$ in a row by a horizontal arc of constant $\theta_2$ and join each $O$ and $X$ in a column by a vertical arc of constant $\theta_1$ mod $\frac{2\pi}{p}$. 

\begin{figure}
\centering
\input{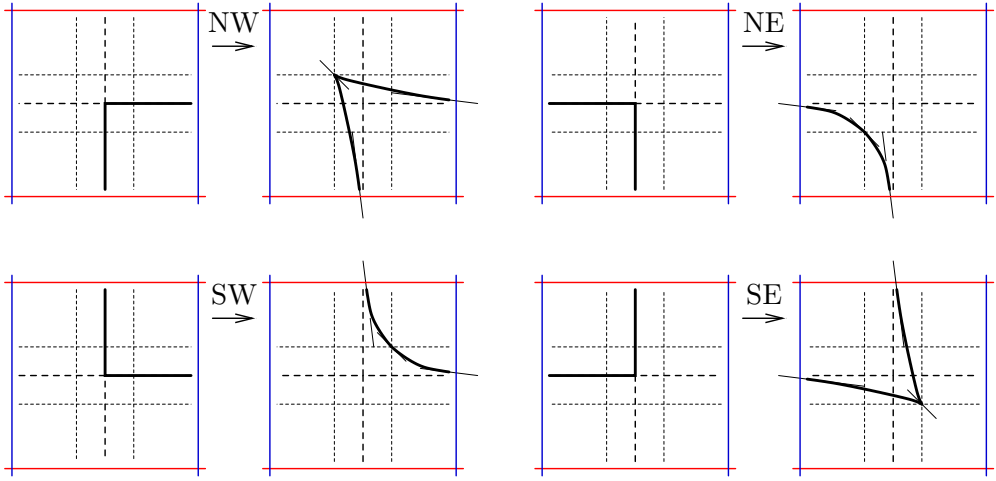}
\caption{The smoothings and cuspings of corners in the rectilinear projection associated to a grid diagram $G_K$ on $\Sigma$.}
\label{fig:smoothing}
\end{figure}

We now perturb this rectilinear projection to yield a toroidal front projection as follows.  The corners of the rectilinear diagram coincide with the $O$'s and $X$'s.  Replace each NW and SE corner with a semi-cubical cusp just outside the corner; replace each NE and SW corner with a rounding just inside the corner.  See Figure~\ref{fig:smoothing}.    These may be done so that they are tangent to the induced line field (of slope $-1$) on $\Sigma$ and so that the associated Legendrian curve intersects $\Sigma$ near the original $O$'s and $X$'s.  Furthermore the smoothed and cusped corners may now be joined by curves with finite negative slopes (in $(-\infty,-1]$ for corners in a column and in $[-1,0)$ for corners in a row) producing a toroidal front.  These choices may be made so that the toroidal front is arbitrarily close to the original rectilinear projection.  Furthermore, the toroidal front isotopy class of the result is unique.  
\begin{figure}
\centering
\input{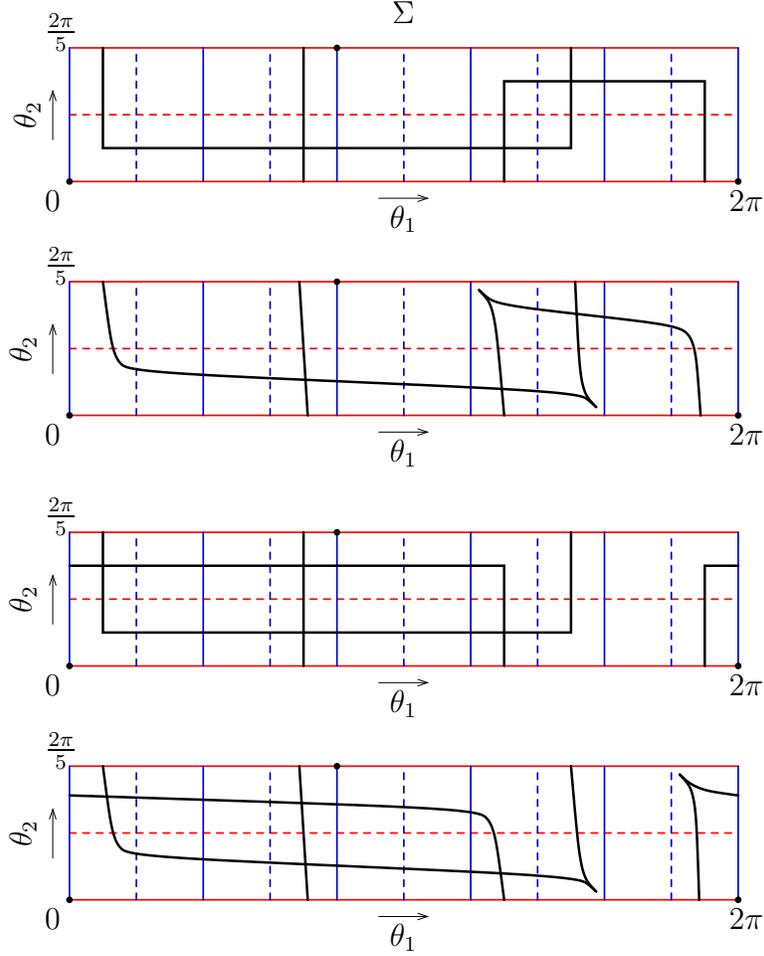}
\caption{Here are two rectilinear projections for the same grid diagram of grid number $2$ for $L(5,2)$ and their associated toroidal front diagrams on $\Sigma$.}
\label{fig:example}
\end{figure}
\end{proof}

By Proposition~\ref{prop:frontgivesleg}, there is a Legendrian link associated to the front obtained above.  Figure~\ref{fig:example} shows two toroidal front diagrams obtained by smoothing two rectilinear projections obtained from the same grid diagram.

\subsection{Legendrian Isotopy Class Invariance of Constructed Toroidal Front}\label{subsec:choices}

\begin{proposition}\label{prop:rectilinearchoice}
The Legendrian isotopy class of the link obtained from a grid diagram $G_K^*$ does not depend upon the choice of rectilinear projection used to define the toroidal front.
\end{proposition}

\begin{figure}
\centering
\input{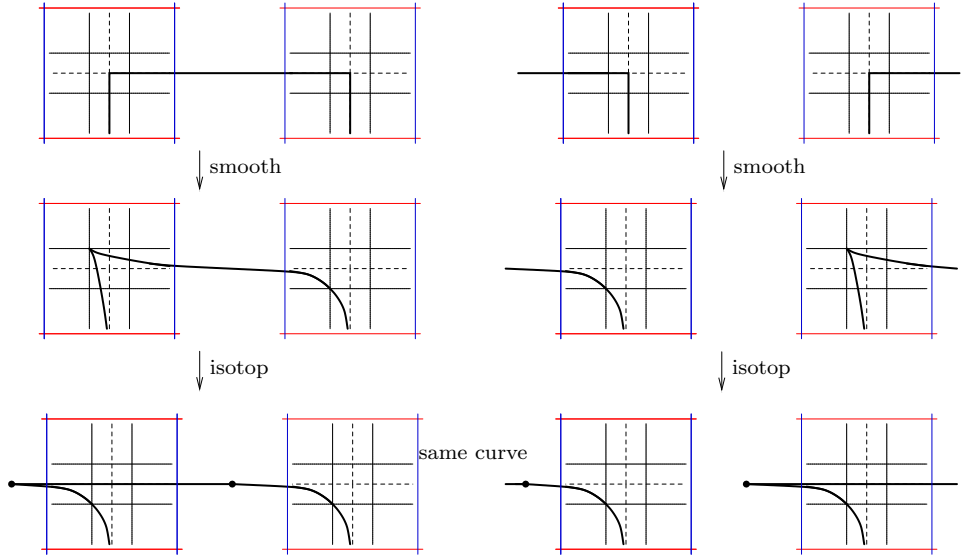}
\caption{ The left- and right-hand sides of this figure show the Legendrian isotopy equivalence of the Legendrian curves obtained from the two choices of horizontal arc in the rectilinear projection connecting an $O$ and an $X$ in a row of a grid diagram.  The rectilinear projections on the upper level may both be smoothed to toroidal fronts as in the middle level.  These two toroidal fronts can be moved via toroidal front isotopies to diagrams (bottom row) which, though not toroidal fronts, represent the same, smooth Legendrian link.}  
\label{fig:isotopydiska}
\end{figure}

\begin{figure}
\centering
\input{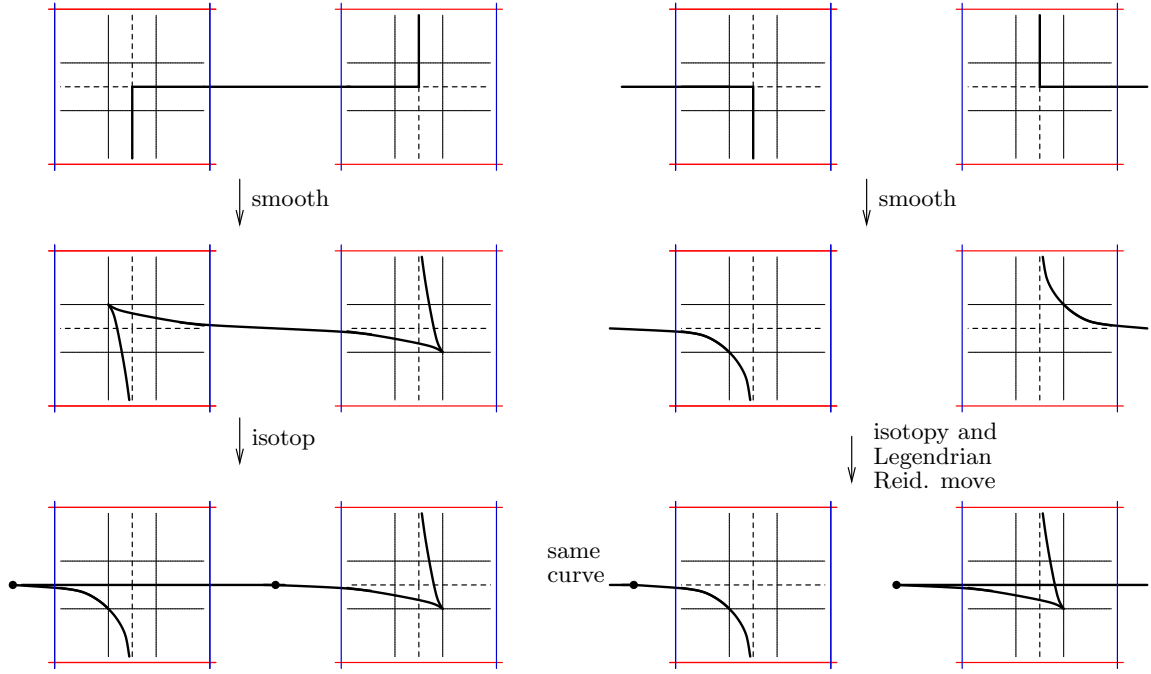}
\caption{ See the caption for Figure~\ref{fig:isotopydiska}.  To pass from the rightmost middle row to the bottom row, one needs to use a local Legendrian Reidemeister I move (see Section \ref{sec:Equiv}) in addition to a toroidal front isotopy.  See Figure~\ref{fig:isotopydiskc} for more details.}
\label{fig:isotopydiskb}
\end{figure}

\begin{figure}
\centering
\input{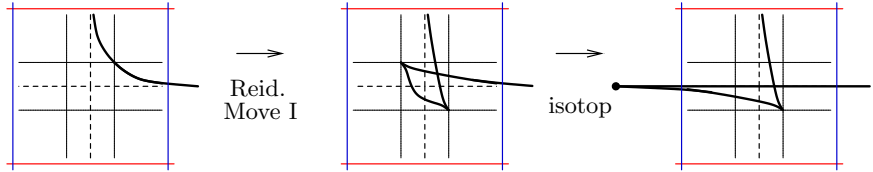}
\caption{  An intermediate isotopy of the lower two right-hand regions of Figure~\ref{fig:isotopydiskb} is through a Legendrian Reidemeister I move.}
\label{fig:isotopydiskc}
\end{figure}

\begin{proof}
It suffices to show that the toroidal fronts obtained from two rectilinear projections differing in a single row represent Legendrian isotopic links; the case for columns is completely analogous.  The top rows of Figures~\ref{fig:isotopydiska} and \ref{fig:isotopydiskb} show the two choices of horizontal arcs for an $O$ and an $X$ in a row in a rectilinear projection.  Figure~\ref{fig:isotopydiska} shows when the vertical arcs are incident to the horizontal arcs from the same side;  Figure~\ref{fig:isotopydiskb} shows when the vertical arcs are incident to the horizontal arcs from opposite sides.  The middle rows illustrate the corresponding toroidal fronts obtained from the rectilinear projections as in Lemma~\ref{lem:griddiagramtotoroidalfront}.

A toroidal front isotopy naturally induces a Legendrian isotopy of the corresponding Legendrian link; we proceed by moving the two toroidal fronts, via toroidal front isotopies, to limiting diagrams which, though not toroidal fronts,  still represent the same smooth Legendrian link.  These two diagrams, pictured in the bottom rows of Figures \ref{fig:isotopydiska} and \ref{fig:isotopydiskb}, both naturally represent a Legendrian link which passes through the core of the solid torus.  The projection has slope $\frac{d\theta_2}{d\theta_1} = 0$ at the two black dots, $p_1, p_2$, which are positioned distance $\theta_1 = \pi$ away from each other on $\Sigma$.  In other words, in both diagrams, the entire horizontal arc connecting $p_1$ and $p_2$ represent a single point, $p$, on the core of the solid torus, and the points near $p_1$ and $p_2$ with nonzero slope are the projections of points nearby $p$ (which project to opposite sides of the Heegaard torus, as pictured).  This Legendrian link will be smooth at all points away from $p$, since the toroidal front projection is smooth at all points away from $p_1$, $p_2$.  To ensure that the Legendrian link is smooth at $p$, one need only arrange for all higher order derivatives of $\theta_2$ with respect to $\theta_1$ to match at $p_1$, $p_2$, which can be done via smooth toroidal front isotopy.
 
Different choices of horizontal or vertical arcs thus produce toroidal fronts of Legendrian isotopic Legendrian links.
\end{proof}

\begin{remark} It is natural to view a rectilinear projection for a Legendrian link as a front diagram for a certain bivalent Legendrian graph whose smoothings produce the toroidal fronts we have been discussing.  Note that radial arcs in $L(p,q)$, i.e., those whose tangent vectors satisfy $\frac{\partial}{\partial \theta_1} = \frac{\partial}{\partial \theta_2} = 0$, are Legendrian.  Therefore, a grid diagram defines a bivalent Legendrian graph made up of a union of radial trajectories between the cores of the two solid tori, passing through $\Sigma$ at the $2n$ basepoints.  A horizontal arc on $\Sigma$ has slope $0$, and thus the points of its interior all define a single point on the core of $V^\alpha$ at radius $r_1=0$, irrelevant $\theta_1$, and specified $\theta_2$.  Similarly a vertical arc has slope $\infty$ and thus the points of its interior all define a single point on the core of $V^\beta$ at radius $r_1=1$, specific $\theta_1$ mod $\frac{2\pi}{p}$, and irrelevant $\theta_2$.  At the endpoints of these horizontal and vertical arcs are the $O$'s and $X$'s.  Since the one-sided tangencies at these points sweep through negative slopes between horizontal and vertical, these points define the radial Legendrian arcs $L_O$ and $L_X$.
\end{remark}

To complete the correspondence between grid diagrams and Legendrian links, we now need only show that any toroidal front diagram is isotopic, through toroidal fronts, to one obtained from a grid diagram:

\begin{lemma}\label{lem:toroidalfronttogriddiagram}
Given a toroidal front diagram on $\Sigma$, there exists a grid diagram $G_K^*$ representing the associated Legendrian link.
\end{lemma}

\begin{proof}
Let $\gamma$ be a toroidal front on $\Sigma$.  Perturb $\gamma$ slightly by a toroidal front isotopy so that all tangencies of slope $-1$ are isolated.  If $\gamma$ has a tangency of slope $-1$ at a non-cusp point with neighboring points having slopes either all strictly greater than $-1$ or  all strictly less than $-1$, then a slight isotopy of $\gamma$ in this neighborhood will eliminate the $-1$ sloped tangency.  Isotop to remove all such tangencies.   Thus on the arcs along $\gamma$ between consecutive tangencies of slope $-1$ and cusps, the tangencies to $\gamma$ have slopes in the range $(-\infty, -1)$ or $(-1,0)$.

Similarly, if $\gamma$ has a cusp with neighboring points having slopes either all strictly greater than $-1$ or  all strictly less than $-1$, then we may arrange, via a slight isotopy, that $\gamma$ instead looks locally like Figure~\ref{fig:turningcusps} near the cusp.  As a consequence, the tangency at each cusp has slope $-1$.  Furthermore, if $p$ is a point of $\gamma$ whose tangent has slope $-1$, then  points on $\gamma$ near $p$ to one side have slopes $> -1$, and points near $p$ to the other side will have slopes $< -1$. 

\begin{figure}
\centering
\input{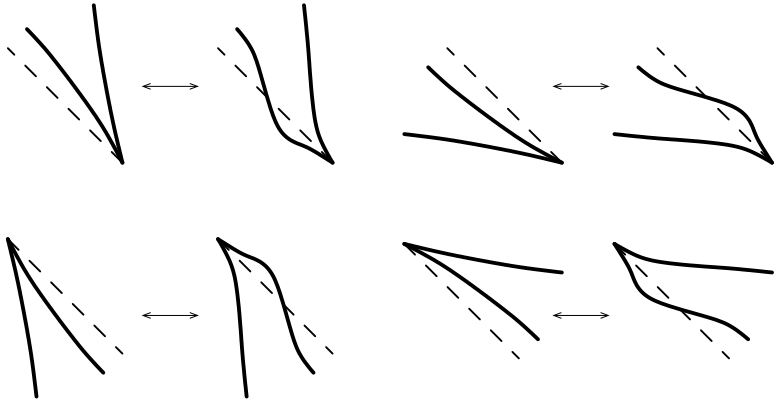}
\caption{Arrange cusps by Legendrian isotopies to have slope $-1$.}
\label{fig:turningcusps}
\end{figure}

Mark each tangency of slope $-1$.  We will refer to consecutive markings as a {\em horizontal pair} if the arc of $\gamma$ joining them has slopes in the range $(-1,0)$; similarly, we will call consecute markings a {\em vertical pair} if the arc joining them has slopes in the range $(-\infty,-1)$.  Arrange, by toroidal front isotopy, that the arc joining each horizontal pair is very close to slope $0$.

It is now straightforward to approximate the toroidal front by a rectilinear projection composed of horizontal and vertical segments in such a way that the associated toroidal front smoothing is isotopic, through toroidal fronts, to the original toroidal front.  Note that each horizontal pair can be approximated by a single horizontal segment, while the vertical pairs may be approximated by zig-zags which, when perturbed to yield a toroidal front diagram, produce no cusps.
\end{proof}

\subsection{Constructing a Planar Front from a Grid Diagram} \label{sec:planarfront}
For $K$ a link in $S^3$, \cite{GT0611841} describes how one associates to a grid diagram $G_K$ for $K$ a standard planar front diagram for a Legendrian representative of the topological mirror of $K$, $m(K)$.  In this section, we prove that the planar front diagram described in \cite{GT0611841} can actually be viewed as a toroidal front diagram supported in a planar region on the Heegaard torus for $S^3$ with the reverse orientation.  As a straightforward consequence of our argument, one can associate to any planar subset of a toroidal front diagram a planar front diagram representing a Legendrian tangle in $(B^3,\xi_{\st})$.  This will be of use in defining Legendrian Reidemeister moves for lens spaces, as we do in Section \ref{sec:Equiv}.

\begin{figure}
\centering
\input{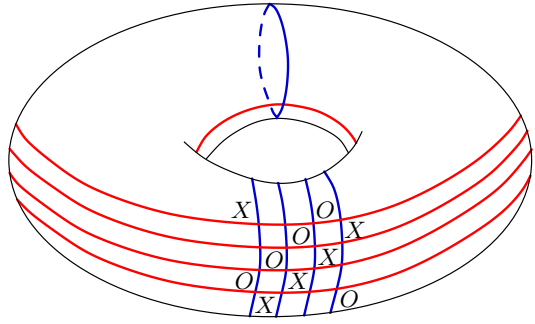}
\caption{A grid diagram whose basepoints are supported in a small planar region of the Heegaard torus.}
\label{fig:RectGrid}
\end{figure}

Let $G_K$ be a grid diagram for $K \subset S^3$ and $G_K^*$ the grid diagram for $m(K) \subset -S^3$.  By distinguishing an $\alpha$ curve and a $\beta$ curve on $G_K^*$, we may choose the horizontal and vertical arcs connecting the $X$'s and $O$'s so that they are disjoint from the two distinguished curves.  (This may be done since each $\alpha$ and $\beta$ curve on a genus $1$ Heegaard diagram for $S^3$ intersect just once.)  We thus obtain a toroidal projection that is confined to a planar subset of the torus.  By the results of the previous section, we know that the Legendrian isotopy class of the knot is unaffected by this choice.  See Figure \ref{fig:RectGrid}.

We may assume our distinguished $\alpha$ and $\beta$ curves correspond to the circles $\theta_2  = \pm \pi$ and $\theta_1 = \pm \pi$ on the Heegaard torus so that all the basepoints of our grid diagram have coordinates $(\theta_1,\theta_2) \in (-\pi,\pi)\times(-\pi,\pi)$.  After smoothing to obtain a toroidal front for the Legendrian link associated to the grid diagram, any tangency is neither horizontal nor vertical. Then using our familiar coordinates $(r_1, \theta_1, \theta_2)$ for points in $S^3$, our Legendrian link is supported in the open tetrahedron 
\[W = \{(r_1, \theta_1, \theta_2)\,\, |\,\, r_1 \in (0,1), \theta_1 \in (-\pi,\pi), \theta_2 \in (-\pi,\pi)\}.\]
The map $f\colon W \to \R^3$ defined by 
\[f(r_1, \theta_1, \theta_2) = \left(\frac{1}{2}(\theta_1 - \theta_2)\, ,\, 1- 2r_1^2\, ,\, \frac{1}{2}(\theta_1 + \theta_2)\right) \]
induces a contactomorphism from ${\xi_{\st}} \vert W$ to ${\xi_{\R^3}} \vert {f(W)}$.  Recall $\xi_{\st} = \ker(r_1^2 d\theta_1 + (1-r_1^2) d\theta_2)$ is the standard contact structure on $S^3$ (and hence on $W$) while $\xi_{\R^3} = \ker (dz - ydx)$ is the standard contact structure on $\R^3$.

The square $W_{r_1=1/\sqrt{2}}$ that is the complement of the distinguished $\alpha$ and $\beta$ curves is mapped to the open diamond in the $xz$--plane of $\R^3$ with vertices at $(\pm \pi, 0, 0)$ and $(0,0, \pm \pi)$.  In this manner the map $f$ carries a toroidal front for a Legendrian link in $S^3$ to a standard front for a Legendrian link in $\R^3$.  Indeed, the image $f(W)$ is the open rectangular solid obtained by sweeping this open diamond in the $xz$--plane along the interval $(-1,1)$ of the $y$--axis.

By Legendrian isotopies, we may arrange that any Legendrian link in $\R^3$ is contained in $f(W)$.  Accordingly, by scaling and vertical compression, we may arrange a front to be contained within this diamond so that every tangent line has slope in the range $(-1,1)$.  Hence (after any necessary scaling and compression) the map $f^{-1}$ carries a front for a Legendrian link in $\R^3$ to a toroidal front.

In other words, $f$ sends our toroidal front to a standard planar front.  Up to planar front isotopy, this is equivalent to rotating the planar subset of the toroidal front $45^{\circ}$  counterclockwise.  To verify that our construction matches the one desribed in \cite{GT0611841}, begin by rotating a planar rectilinear projection for $G_K$ $90^{\circ}$ clockwise to produce a planar rectilinear projection for $G_K^*$.  This has the effect of changing all crossings of the associated link.  Then rotate the toroidal front $45^{\circ}$ counterclockwise to produce a planar front.  This is precisely the procedure described in \cite{GT0611841}.

\section{Topological and Legendrian Equivalence Under Grid Moves} \label{sec:Equiv}
Now that we have established a relationship between a grid diagram $G_K$ for a link $K \subset L(p,q')$ and a Legendrian representative $\mathcal{L}_K$ of $K$ in $(L(p,q),\xi_{\UT})$, we turn to the question of the uniqueness of this representation.  Our goal is to provide an elementary set of moves, as in \cite{MR1339757},\cite{GT0610559},\cite{GT0611841}, allowing one to move between any two grid diagram representatives of the same topological or Legendrian knot type.

The main theorem of this section is the following:

\begin{theorem}  Let $G_K$ and $G_K'$ be grid diagrams representing smooth links $K$ and $K'$ in $L(p,q')$ (resp., Legendrian links $\mathcal{L}_K$ and $\mathcal{L}_{K'}$ in $(L(p,q),\xi_{\UT})$).  

\begin{enumerate} \label{thm:GridMoves}
\item $K$ and $K'$ are smoothly isotopic iff there exists a sequence of elementary topological grid moves connecting $G_K$ to $G_K'$.
\item $\mathcal{L}_K$ and $\mathcal{L}_{K'}$ are Legendrian isotopic iff there exists a sequence of elementary Legendrian grid moves connecting $G_K$ to $G_K'$.\end{enumerate}
\end{theorem}

\begin{figure}
\centering
\input{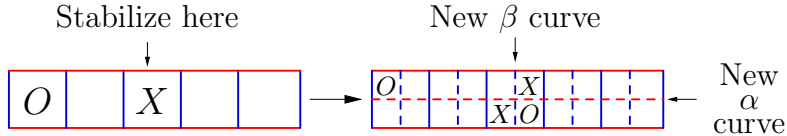}
\caption{An illustration of a NW stabilization at an $X$ basepoint for a grid number 1 knot in $L(5,1)$.  One introduces new $\alpha$-$\beta$ and $X$-$O$ pairs in the rectangle specified by the chosen basepoint.A destabilization is the inverse of this operation.}
\label{fig:Stab}
\end{figure}

\begin{figure}
\centering
\input{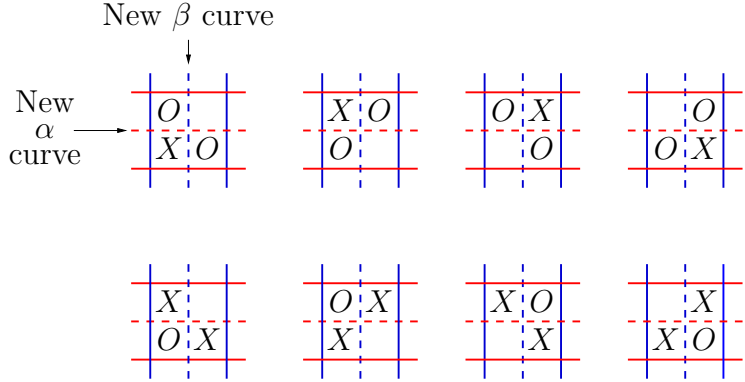}
\caption{Local pictures for all $8$ different types of stabilizations.  Starting from the left on the top row, these are denoted O:NE, O:SE, O:SW, and O:NW.  Along the bottom row, we have X:NE, X:SE, X:SW, and X:NW.}
\label{fig:AllStabs}
\end{figure}

\begin{figure}
\centering
\input{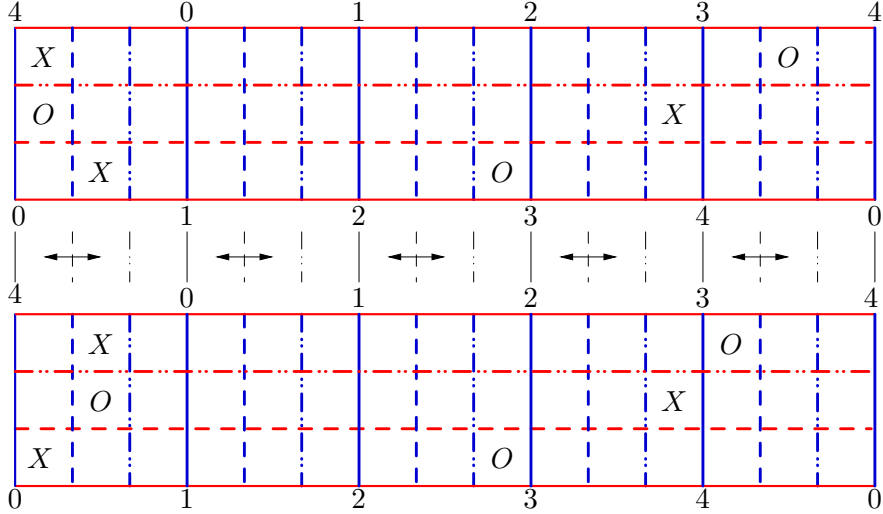}
\caption{A before (top) and after (bottom) snapshot of a first and second column commutation on a grid number 3 diagram for a link in $L(5,1)$.  Since the basepoint pairs in these columns do not interleave, we can exchange the markings as shown.  A row commutation is the obvious analogue for two adjacent rows.  In general, a row or column commutation can be performed on any two adjacent rows or columns as long as the markings in the two columns do not interleave.}
\label{fig:Comm}
\end{figure}

Figures \ref{fig:Stab}, \ref{fig:AllStabs}, and \ref{fig:Comm} illustrate the elementary topological grid moves: (de)stabilizations and commutations.  These moves look locally like those described in \cite{GT0611841} for $S^3$ knots.  The elementary Legendrian grid moves form a subset of the elementary grid moves, including commutations and (de)stabilizations of types X:NW, X:SE, O:NW, O:SE.  We will say that two grid diagrams are topologically (resp., Legendrian) grid-equivalent if they are related by a finite sequence of elementary topological (resp., Legendrian) grid moves.

\begin{figure}
\centering
\input{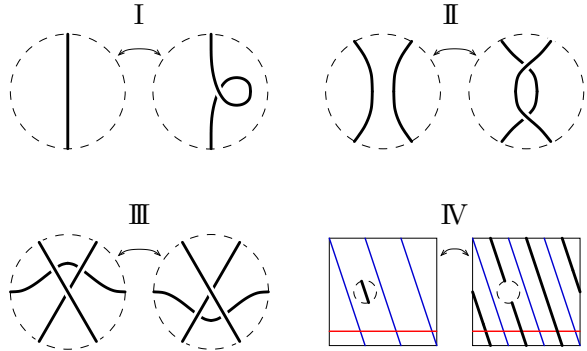}
\caption{An illustration of the lens space Reidemeister moves.  Moves I, I\!I, I\!I\!I look locally like the $S^3$ Reidemeister moves, but should be thought of here as occurring in a sufficiently small disk on the Heegaard torus.  Move I\!V corresponds to replacing an arc of the boundary of a meridional disk (in either solid handlebody) with the complementary arc.  Here, we have pictured this move for a meridional disk in $V^{\beta}$ in the lens space $L(3,1)$. }
\label{fig:Reid}
\end{figure}

\begin{proof} 
{\flushleft {\bf Item 1: Smooth Isotopy}}
The Reidemeister theorem for links in $S^3$ states that two links are smoothly isotopic iff their planar projections are related by planar isotopies and Reidemeister moves.  Since each such move corresponds to a local isotopy (supported in  $\R^3$), it is immediate that grid-equivalent grid diagrams for links in lens spaces represent smoothly isotopic links.  Furthermore, the observation that all Reidemeister moves are local, coupled with the natural correspondence between planar subsets of toroidal projections of lens space links and standard planar projections of tangles in $\R^3$, immediately implies a version of Reidemeister's theorem for lens space links:

\begin{proposition} Two smooth links $K_1$ and $K_2$ in $L(p,q)$ are smoothly isotopic iff their toroidal projections are related by a sequence of smooth isotopies and Reidemeister moves $I, I\!I, I\!I\!I,$ and $IV$, picture in Figure \ref{fig:Reid}.
\end{proposition}

To prove that any two smoothly isotopic links $K_1$ and $K_2$ have topologically grid-equivalent grid diagrams, we use the argument outlined by Dynnikov in \cite{MR2232855}.  First, we represent $K_1$ and $K_2$ by (rectilinear) toroidal projections, using Lemma 4.2 of \cite{GT07100359}.  By the Reidemeister theorem for lens space links, we can move from one toroidal projection to the other by a sequence of smooth isotopies and moves of type I - I\!V.  By the same arguments used in the proofs of Lemma 4.2 and Proposition 4.3 of \cite{GT07100359}, we can approximate each stage of this process using a grid diagram.  Provided that each intermediate step is sufficiently simple (subdivide the compact isotopy further if not), it is easy to verify that each step can be accomplished using elementary grid moves.  Reidemeister move I\!V does not require a grid move; rather, one chooses an alternate projection of the link to the Heegaard torus.  Figures \ref{fig:ReidI}, \ref{fig:ReidII}, and \ref{fig:ReidIII} enumerate all possible versions of the other Reidemeister moves, indicating how they can be obtained via grid moves.  

\begin{figure}
\centering
\input{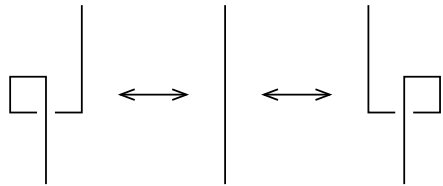}
\caption{Reidemeister I via grid moves.}
\label{fig:ReidI}
\end{figure}

\begin{figure}
\centering
\input{ReidII.pstex_t}
\caption{Reidemeister I\!I via grid moves.}
\label{fig:ReidII}
\end{figure}

\begin{figure}
\centering
\input{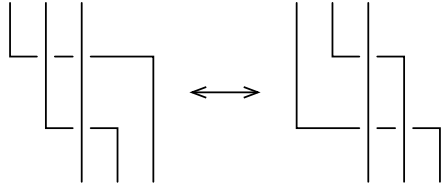}
\caption{Reidemeister I\!I\!I via grid moves.}
\label{fig:ReidIII}
\end{figure}

{\bf Item 2: Legendrian Isotopy}
In \cite{MR1186009}, Swiatkowski proves that two Legendrian links in $S^3$ are Legendrian isotopic iff their corresponding planar fronts are related by a sequence of planar isotopies and Legendrian Reidemeister moves as in Figure~\ref{fig:LegReid}.  Since each of his Reidemeister moves corresponds to a local Legendrian isotopy (i.e., takes place in a Darboux ball), his result will immediately imply that Legendrian grid-equivalent grid diagrams represent Legendrian-isotopic lens space links, once we understand the relationship between planar subsets of toroidal fronts and standard planar fronts.  This relationship, explained in detail in Section \ref{sec:planarfront}, coupled with Swiatkowski's result, yields a Legendrian Reidemeister theorem for Legendrian links in $(L(p,q),\xi_{\UT})$:

\begin{figure}
\centering
\input{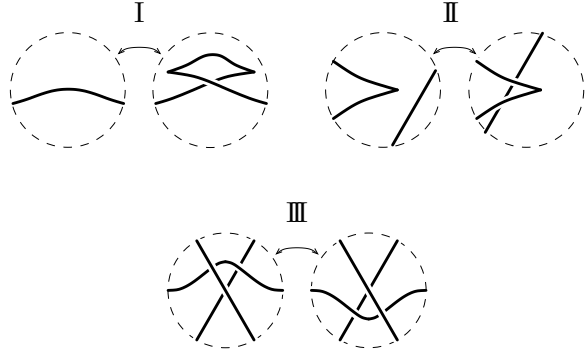}
\caption{Legendrian Reidemeister moves for planar fronts.}
\label{fig:LegReid}
\end{figure}

\begin{figure}
\centering
\input{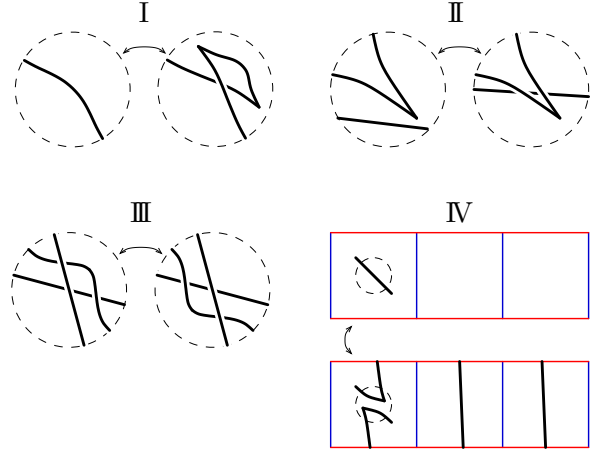}
\caption{Legendrian Reidemeister moves for toroidal fronts.}
\label{fig:LegReidtoroidal}
\end{figure}

\begin{proposition} Two toroidal front diagrams represent Legendrian-isotopic links in $(L(p,q),\xi_{\UT})$ iff they can be connected by a sequence of toroidal front isotopies, Legendrian Reidemeister moves, and Legendrian slides across meridional disks (which we will call Legendrian Reidemeister moves of type $IV$), as in Proposition \ref{prop:rectilinearchoice}.
\end{proposition}

The Legendrian Reidemeister moves I - I\!V on a toroidal front are illustrated in Figure~\ref{fig:LegReidtoroidal}.

To prove that Legendrian-isotopic links $K_1$ and $K_2$ have Legendrian grid-equivalent grid diagrams, one first uses Proposition \ref{prop:leggivesfront} to represent both by toroidal fronts on $\Sigma$.  By the Legendrian Reidemeister theorem for links in lens spaces, we can move between the two fronts by a sequence of toroidal front isotopies and Legendrian Reidemeister moves.  An argument exactly as in the proof of Lemma \ref{lem:toroidalfronttogriddiagram} produces a dual grid diagram $G_{K_{i}}^*$ for each stage of the Legendrian isotopy.  It is then straightforward, using arguments analogous to those given in the verification of topological grid-equivalence, to show that the associated grid diagrams $G_{K_i}$ ({\em not} the dual grid diagrams $G_{K_i}^*$) associated to each stage of the Legendrian isotopy are Legendrian grid-equivalent.
\end{proof}

\section{The Thurston-Bennequin invariant and bounds on grid number} \label{sec:ClassInv}
Let $\mathcal{L}_K$ be an oriented Legendrian link in a contact rational homology sphere $(Y,\xi)$ repesenting the topological oriented link $K$ with components $K_1, \dots, K_\ell$.  We begin by defining the Thurston-Bennequin number, $\tb(\mathcal{L}_K)$, a classical Legendrian invariant.  Traditionally $\tb$ is defined when $K$ is a null-homologous link in a contact $3$--manifold, see e.g.\ \cite{MR1445553}.  The definition extends in a natural way to rationally null-homologous oriented links. 

\subsection{Contact framing, Seifert framing, and the Thurston-Bennequin number}

Let $N(K)$ denote a small tubular neighborhood of the link $K$ with $N(K_i)$ denoting the component that is a neighborhood of $K_i$.  Let $\mu_i$ be an oriented meridian of the closure of $N(K_i)$ so that it links $K_i$ once positively.  Typically a framing for $K_i$ is a choice of a slope on $\bdry N(K_i)$ that algebraically intersects the meridian $\mu_i$ once.  We shall relax this notion of framing so that it may be a collection of parallel slopes on $\bdry N(K_i)$ that each algebraically intersect the meridian $\mu_i$ a fixed number of times.

The Thurston-Bennequin number $\tb(\mathcal{L}_K)$ measures the discrepancy between the {\em contact framing} and the {\em Seifert framing} of the oriented Legendrian link $\mathcal{L}_K$.   

\begin{definition} 
The {\em contact framing} for $\mathcal{L}_K$ is the tuple 
\[{\boldsymbol \gamma} = ([\gamma_1], \ldots, [\gamma_\ell]) \in \bigoplus_{i=1}^\ell H_1(T^2_i;\Z)\]
defined by pushing the oriented components $K_i$ of $K$ into $T^2_i = \bdry(Y-N(K_i)) \subset \bdry(Y-N(K))$ along the contact planes.
\end{definition}

\begin{definition}
The {\em Seifert framing} for $K$ is the tuple
\[{\boldsymbol \lambda} = ([\lambda_1], \ldots, [\lambda_\ell]) \in \bigoplus_{i=1}^\ell H_1(T^2_i;\Z)\]
such that 
\begin{enumerate}
\item $\lambda_i$ is a collection of parallel, coherently oriented, simple closed curves on $T^2_i=\bdry N(K_i)$,
\item $\mu_i \cdot \lambda_i = d$ for each $i$ where $d$ is the order of $[K] \in H_1(Y;\Z)$, and
\item via the induced maps on homology coming from the inclusions $T_i^2 \hookrightarrow Y-N(K)$, we have $\sum_{i=1}^{\ell} \lambda_i = 0 \in H_1(Y-N(K);\Z)$.
\end{enumerate}
\end{definition}

It is clear that the contact framing is well-defined, but it may not be so clear that  the Seifert framing is well-defined for non-nullhomologous oriented links in rational homology spheres.

\begin{lemma}
The Seifert framing is well-defined.
\end{lemma}

\begin{proof}
Given an oriented link $K = K_1 \cup \dots \cup K_\ell$ in a rational homology sphere $Y$ with oriented meridians $\mu_i$ as above, observe that $H_1(Y-N(K);\Q) \cong \Q^\ell$.  Moreover the homology classes of the oriented meridians form a basis $\{ [\mu_1], \dots, [\mu_\ell]\}$.  Set 
\[{\boldsymbol \mu} = ([\mu_1], \dots, [\mu_\ell]) \in \bigoplus_{i=1}^{\ell} H_1(\bdry (Y-N(K));\Q).\]  

We construct the Seifert framing as follows.  For each component $K_i$ of $K$ pick an oriented push-off,  $\lambda_i'$, of $K_i$ on $\bdry N(K_i)$.  Note that $\lambda_i'$ is a framing in the usual sense.  Write 
\[{\boldsymbol \lambda'} = ([\lambda_1'], \dots, [\lambda_\ell']) \in \bigoplus_{i=1}^{\ell} H_1(\bdry(Y-N(K));\Q).\]

Since each $[\lambda_i]$ may be expressed as a $\Q$--linear combination of the $[\mu_j]$, there exists an $\ell \times \ell$ matrix $A$ with rational coefficients such that
\[ {\boldsymbol \lambda'} = A {\boldsymbol \mu}.\]

Let $d$ be the smallest positive integer such that $dA$ has integral coefficients, and let ${\boldsymbol \lambda}^{\Sigma}$ denote the image of ${\boldsymbol \lambda}$ under the natural map \[ \bigoplus_{i=1}^{\ell} H_1(Y-N(K);\Z) \rightarrow H_1(Y-N(K);\Z)\] which sends \[ (\lambda_1, \ldots, \lambda_{\ell}) \rightarrow \sum_{i=1}^\ell \lambda_i.\]

Then, we see that \[d{\boldsymbol \lambda'}^\Sigma - d(A{\boldsymbol \mu})^\Sigma = 0 \in H_1(Y-N(K);\Z).\]

Furthermore, by collecting like terms, one produces a unique element \[ {\boldsymbol \lambda} = (\lambda_1, \ldots, \lambda_{\ell}) \in \bigoplus_{i=1}^\ell H_1(T_i^2;\Z). \]  More precisely, \[ \lambda_j = \lambda_j' - \sum_{i=1}^{\ell} A_{ij} \mu_j. \]

To see that the resulting element of $\bigoplus_{i=1}^{\ell} H_1(T_i^2;\Z)$ is independent of the original choices of the push-offs $\lambda_i'$, simply note that for each $K_i$, any other choice, $\lambda_i''$, of push-off will differ from $\lambda_i'$ by an integer multiple of $\mu_i$.  Since the $[\mu_i]$ form a basis for $H_1(Y-N(K);\Q)$, this has a predictable compensatory effect on the corresponding matrix $A''$ expressing ${\boldsymbol \lambda''}$ in terms of the $\mu_i$; hence, \[{\boldsymbol \lambda''} - A''{\boldsymbol \mu} = {\boldsymbol \lambda'} - A{\boldsymbol \mu},\] as desired.

Observe that $d$ is the order of $[K]$ in $H_1(Y;\Z)$.
\end{proof}

\begin{definition}
Let ${\boldsymbol \lambda} \in \bigoplus_{i=1}^{\ell} H_1(T_i^2;\Z)$ be a Seifert framing, as constructed above.  Then ${\boldsymbol \lambda}^\Sigma$ bounds an oriented surface, $F$, properly embedded in $Y-N(K)$.  We will call any such surface, $F$, a {\em (rational) Seifert surface} for the link $K$.  One may also consider contracting $N(K)$ radially to $K$ so that the interior of $F$ is embedded while $\bdry F$ is a $d$--fold cover of $K$.  (Again, $d$ is the order of $[K]$ in $H_1(Y;\Z)$.)
\end{definition}

\begin{remark} If $Y$ is not a rational homology sphere, then there may be a $\Q$--linear dependence among the homology classes of the meridians of $K$ in $H_1(Y-N(K);\Q)$.  Such an occurrence would cause the Seifert framing to be ill-defined.  As an example, consider the null homologous link $S^1 \times {x} \cup S^1 \times {y}$ in $S^1 \times S^2$ for distinct points $x,y \in S^2$. The ambiguity may be resolved by specifying the $2$nd homology class of the Seifert surface.

\end{remark}

We are now ready to define the Thurston-Bennequin number of a Legendrian link $\mathcal{L}_K$.

\begin{definition} \label{defn:tb}
Let $\mathcal{L}_K$ be a Legendrian representative of a link $K = \bigcup_{i=1}^\ell K_i$ of order $d$ in a contact rational homology sphere $(Y,\xi)$, ${\boldsymbol \gamma}$ its contact framing, and ${\boldsymbol \lambda}$ its Seifert framing. Then the {\em Thurston-Bennequin number} of $\mathcal{L}_K$, $\tb(\mathcal{L}_K) \in \Q$, is:
\[\tb(\mathcal{L}_K) = \frac{1}{d} ({\boldsymbol \gamma} \cdot {\boldsymbol \lambda}).\]
\end{definition}

We have abused notation slightly in the above definition, thinking of ${\boldsymbol \gamma} \cdot {\boldsymbol \lambda}$ as an inner product of algebraic intersection numbers.  More precisely,
\[{\boldsymbol \gamma} \cdot {\boldsymbol \lambda} := \sum_{i=1}^\ell [\gamma_i] \cdot [\lambda_i],\]
where ``$\cdot$'' on the right hand side refers to algebraic intersection number.

Note that an overall reversal of orientations on $\mathcal{L}_K$ will not alter $\tb(\mathcal{L}_K)$.  However changing the orientation of one component of a multi-component link could drastically alter the Seifert framing and thus change the Thurston-Bennequin number.

The following relationship between $\tb(\mathcal{L}_K)$ and $\tb(\mathcal{L}_{\widetilde{K}})$ under a contact cover of degree $m$ is immediate.

\begin{lemma}\label{lemma:LiftTBR} 
Let $\mathcal{L}_K$ be a Legendrian link in the contact rational homology sphere $(Y,\xi)$ and $\mathcal{L}_{\widetilde{K}}$ its lift to a Legendrian link in an $m$--fold  contact cover $(\widetilde{Y},\widetilde{\xi})$.  Then 
\[\tb(\mathcal{L}_K) = \frac{1}{m} \tb(\mathcal{L}_{\widetilde{K}}).\]
\end{lemma}

\begin{proof}
Note that $\widetilde{Y}$ will also be a rational homology sphere.
By construction, ${\boldsymbol \lambda}$ is the image, under the map induced by the projection $\pi\colon \widetilde{Y} \to Y$, of the Seifert framing $\widetilde{\boldsymbol \lambda}$ for $\widetilde{K}$.  Similarly, ${\boldsymbol \gamma}$ is the image, under $\pi_*$, of the contact framing $\widetilde{\boldsymbol {\gamma}}$ for $\widetilde{K}$. 

 Since
\[ {\boldsymbol \gamma} \cdot {\boldsymbol \lambda} = (\pi_* \widetilde{{\boldsymbol \gamma}}) \cdot (\pi_* \widetilde{\boldsymbol {\lambda}}) = \pi_*(\widetilde{\boldsymbol \gamma}\cdot \widetilde{\boldsymbol \lambda}),\]
 we see that 
 \[ \tb(\mathcal{L}_K) =\frac{1}{m} \tb(\mathcal{L}_{\widetilde{K}})\]
  as desired.
\end{proof}

\subsection{Computing $\tb(\mathcal{L}_K)$ from Grid diagrams}
Since $\tb$ behaves simply when taking contact covers (Lemma \ref{lemma:LiftTBR}), we can compute $\tb(\mathcal{L}_K)$ from combinatorial data on $G_K$ for any grid diagram representing $\mathcal{L}_K$.  Recall that $G_K$ represents $\mathcal{L}_K$ if a toroidal front diagram for $\mathcal{L}_K$ can be obtained from $G_K$ by the procedure described in the proof of Theorem \ref{thm:Legrealize}.  Although this combinatorial data is most easily expressed in the language of Heegaard Floer homology, we will now mention Heegaard Floer homology only to the extent necessary to identify the relevant combinatorial data.  See \cite{GT07100359} for more details.

Associated to a grid number $n$ grid diagram $G_K$ for $K \in L(p,q)$ is a combinatorial Heegaard Floer chain complex, $CF^-(G_K)$, whose generators correspond to one-to-one matchings between the $\alpha$ and $\beta$ curves.  Each generator, ${\bf x}$, furthermore comes equipped with a combinatorially-defined homological grading, ${\bf M}({\bf x}) \in \mathbb{Q}$, as described in Section 2.2.2 of \cite{GT07100359}.  

Now let $z^-$ be the generator of $CF^-(G_K)$ in the SW (lower left) corner of the $X$ basepoints and $z^+$ be the generator of $CF^-(G_K)$ in the NE (upper right) corner of the $X$ basepoints.  Abusing notation so that $G_K$ refers to the Legendrian link represented by $G_K$, we have:

\begin{proposition} \label{prop:tbrformula}
\[
\tb(G_K) = \frac{{\bf M}(z^-) + {\bf M}(z^+)}{2} - d(p,q,q-1) - 1
\]
\end{proposition}

Here, $d(p,q,q-1)$ denotes the correction term $d(-L(p,q),q-1)$ as defined inductively in \cite{MR1957829}.

\begin{proof}
We combine 
Theorem 1.1 from \cite{GT0611841} and Equation 2 from \cite{GT07100359} to conclude that
\begin{eqnarray*}
 {\bf M}(z^-) + {\bf M}(z^+) &=& \frac{1}{p}\left({\bf M}(\widetilde{z}^-) + {\bf M}(\widetilde{z}^+)\right) + 2d(p,q,q-1) + \frac{2(p-1)}{p}\\
  &=& \frac{2}{p}\left(\tb(\widetilde{G}_{K}) + 1 + p\cdot d(p,q,q-1) + (p-1)\right). 
\end{eqnarray*}
Lemma~\ref{lemma:LiftTBR} then implies that 
\[\tb(G_K) = \frac{\tb(\widetilde{G}_{K})}{p} = \frac{{\bf M}(z^-) + {\bf M}(z^+)}{2} - d(p,q,q-1) - 1.\]
\end{proof}

The classical Thurston-Bennequin invariant is particularly relevant to the Berge conjecture because of the following relationship, a generalization of a result of Matsuda \cite{MR2240685}.  Recall that $G_K^*$ denotes the grid diagram dual to $G_K$ (see Section \ref{sec:orientation}).  In particular, if $G_K$ represents a Legendrian link with respect to $(L(p,q),\xi_{UT})$, then $G_K^*$ corresponds to a Legendrian link with respect to $(L(p,q'), \xi_{\UT})$, where $qq' \equiv -1 \mod p$.

\begin{proposition} \label{prop:Matsudabound}
\[\tb(G_K) + \tb(G_K^*) = -\GN(G_K).\]
\end{proposition}

\begin{proof} Let $n = \GN(G_K) = \GN(G_K^*).$  Then $G_K$ and $G_K^*$ lift to grid diagrams for $\widetilde{K} \subset S^3$ and $m(\widetilde{K}) \subset S^3$, respectively, where $m(\widetilde{K})$ denotes the topological mirror of $\widetilde{K}$ in $S^3$.  Both grid diagrams have grid number $np$.  Then, as proved in Section 2 of \cite{MR2240685} (Equation 1 of \cite{GT0612356}), we see that 
\[\tb(\widetilde{G}_K) + \tb((\widetilde{G}_{K})^*) = -np.\]
Coupled with Lemma \ref{lemma:LiftTBR}, this implies that 
\[\tb(G_K) + \tb(G_K^*) = -n,\]
as desired.
\end{proof}

Let $\overline{\tb}(K)$ denote the maximum $\tb$ of a Legendrian representative of $K$.  Recalling that if $G_K$ represents $K \subset L(p,q')$, then $G_K^*$ represents $m(K) \subset L(p,q)$ ($qq' \equiv -1 \mod p$), we see that Proposition \ref{prop:Matsudabound} implies:

\begin{corollary}\label{cor:gridnumber}
\[\GN(K) \geq -\overline{\tb}(K) - \overline{\tb}(m(K))\]
for each link $K \subset L(p,q')$.
\end{corollary}

\section{Calculations}
We conclude with some simple observations about Legendrian realizations of knots in $L(p,q)$ for which some integral surgery yields $S^3$.

\begin{proposition} \label{prop:tbmodp} Let $K$ be a knot in $L(p,q)$ upon which integral surgery yields an integer homology sphere $S$, $K'$ the induced knot in $S$, and $\mathcal{L}_{K}$ a Legendrian representative of $K$ with respect to $(L(p,q),\xi_{\UT})$.  Then $$p \cdot \tb(\mathcal{L}_K) \equiv {\pm 1} \mod p.$$
\end{proposition}

\begin{proof} Let $T^2 = \partial(L(p,q) - \mathcal{L}_K)$ be the torus boundary of a neighborhood of $\mathcal{L}_K$.  We denote by 

\begin{itemize}
\item $\mu \subset T^2$ an oriented meridian of $K$ in $L(p,q)$, 
\item $\gamma \subset T^2$ the contact push-off of $K$,
\item $\mu' \subset T^2$ an oriented meridian of $K'$ in $S^3$,
\item $\lambda' \subset T^2$ an oriented simple closed curve representing the Seifert framing of $K'$ in $S^3$ and satisfying $\mu' \cdot \lambda' = 1$.
\end{itemize}

Let $|K|$ denote the order of $K$ as an element of $H_1(L(p,q);\Z)$.  By definition, 
\[\tb(\mathcal{L}_K) := \frac{1}{|K|}(\gamma \cdot \lambda'),\]
and, hence, since $|K| = p$ for all knots admitting integer homology sphere surgery, 
\[p \cdot \tb(\mathcal{L}_K) = \gamma \cdot \lambda'.\]

By assumption, $\mu = p\mu' \pm \lambda'$ (as elements of $H_1(T^2;\Z)$).  Using the fact that $\mu \cdot \gamma = 1$, we conclude that 
\[\gamma = [kp \mp 1]\mu' \pm k \lambda'\]
for some $k \in \Z$.

But this implies that 
\[p \cdot \tb(\mathcal{L}_K) = kp \mp 1,\]
 hence 
 \[p \cdot \tb(\mathcal{L}_K) \equiv \pm 1 \mod p,\]
 as desired.
\end{proof}

\begin{lemma} \label{lem:contsurgery} Let $K$ be a knot in $L(p,q)$ upon which integral surgery yields an integer homology sphere $S$.  If $\mathcal{L}_K$ is a Legendrian representative of $K$ with respect to $(L(p,q),\xi_{\UT})$, and $p \cdot \tb(\mathcal{L}_K) = kp \pm 1$, then the corresponding contact surgery coefficient yielding $S$ is $\pm k$.
\end{lemma}

\begin{proof}
Let $\mu, \mu', \lambda, \lambda', \gamma$ be as above.

As in the proof of Proposition \ref{prop:tbmodp}, if $\mathcal{L}_K$ satisfies $p \cdot \tb(\mathcal{L}_K) = kp \mp 1$, then we can write:
\begin{eqnarray*}
\gamma &=& (kp \mp 1)\mu' \pm k\lambda'.
\end{eqnarray*}
Since, by assumption, $\mu = p\mu' \pm \lambda'$, we see:
\begin{eqnarray*}
\gamma &=& (kp \mp 1)\mu' \pm k[\pm \mu \mp p\mu')]\\
       &=& \mp \mu' + k\mu.
\end{eqnarray*}
Hence, $\pm \mu' = - \gamma + k\mu,$ implying that the contact surgery coefficient is $- k$.
\end{proof}

Note that if $K \subset L(p,q)$ has an integral surgery yielding $S^3$, then $K \subset L(p,q')$ (the topological mirror of $K$) also has an integral surgery yielding $S^3$.  The following is an easy corollary of the preceding two statements:

\begin{corollary}  \label{cor:zeroorone} If $G_K$ is a $\GN1$ diagram for a Legendrian knot $\mathcal{L}_K \subset (L(p,q'),\xi_{\UT})$ whose underlying topological knot has an integral surgery slope yielding $S^3$  and $G_K^*$ is the dual $\GN1$ diagram representing a Legendrian knot $\mathcal{L}_{K^*}$ in $(L(p,q),\xi_{\UT})$, the corresponding contact surgeries on $\mathcal{L}_K$ and $\mathcal{L}_{K^*}$ which yield $S^3$ are $0$ and $+1$, respectively (or vice versa).
\end{corollary}

\begin{proof}  By Proposition \ref{prop:Matsudabound}, 
\[ \tag{\tagGNOneMatsuda} 
\tb(\mathcal{L}_K) + \tb(\mathcal{L}_{K^*}) = -1.
\]
Exchanging the roles of $\mathcal{L}_K$ and $\mathcal{L}_{K^*}$ if necessary, we know, by Proposition \ref{prop:tbmodp}, that 
\begin{eqnarray*}
\tb(\mathcal{L}_K) &=& \frac{k_1p + 1}{p}, \,\,\mbox{ and}\\
\tb(\mathcal{L}_K) &=& \frac{k_2p -1}{p}, 
\end{eqnarray*}
where $k_1, k_2 \in \Z$.  By Equation  
(\tagGNOneMatsuda), $k_1 + k_2 = -1$.

Since any Stein filling of $S^3$ is diffeomorphic to $B^4$ \cite{MR1171908}, $(L(p,q), \xi_{\UT})$ is Stein fillable for all $L(p,q)$ \cite{MR1786111}, and any negative contact surgery can be turned into a sequence of $-1$ contact surgeries \cite{MR2055048}, we can conclude that $k_1, k_2 \leq 0$, since positive $k_i$ would yield a negative contact surgery coefficient by Lemma \ref{lem:contsurgery}, leading to a non-standard Stein filling of $S^3$.

Therefore, $k_1 + k_2 = -1$ implies that $k_1 = 0$ and $k_2 = -1$, or vice versa.  In other words, contact $0$ surgery on one of the two Legendrian knots and contact $+1$ surgery on the other yields $S^3$.
\end{proof}

\bibliography{Legendrian}
\end{document}